\documentclass[reqno]{amsart}
\usepackage{amssymb}
\usepackage[all]{xy}

\newif\ifpdf\ifx\pdfoutput\undefined\pdffalse\else\pdfoutput=1\pdftrue\fi
\ifpdf
    \usepackage{hyperref}
\else
    \usepackage[hypertex=true]{hyperref}
\fi

\newtheorem{thm}{Theorem}[section]
\newtheorem{lem}[thm]{Lemma}
\newtheorem{cor}[thm]{Corollary}
\newtheorem{prop}[thm]{Proposition}

\newtheorem{introthm}{Theorem}

\theoremstyle{definition}
\newtheorem{rem}[thm]{Remark}
\newtheorem{defn}[thm]{Definition}

\renewcommand{\theenumi}{(\roman{enumi})}

\newfont{\cyrr}{wncyr10}
\def\Sh{\mbox{\cyrr Sh}}

\def\Z{\mathbf{Z}}
\def\Q{\mathbf{Q}}
\def\F{\mathbf{F}}

\def\Zp{\Z_p}
\def\Qp{\Q_p}
\def\Fp{\F_p}

\def\A{\mathcal{A}}

\def\E{\mathcal{E}}
\def\cF{\mathcal{F}}
\def\cG{\mathcal{G}}
\def\cR{\mathcal{R}}
\def\cS{\mathcal{S}}

\def\K{\mathcal{K}}
\def\L{\mathcal{L}}
\def\M{\mathcal{M}}
\def\I{\mathcal{I}}
\def\J{\mathcal{J}}
\def\N{\mathcal{N}}

\def\T{W}

\def\ld{\mathcal{h}}
\def\rd{\mathcal{i}}

\def\m{\mathfrak{m}}
\def\p{\mathfrak{p}}

\def\P{\mathfrak{p}}

\def\Hom{\mathrm{Hom}}
\def\Gal{\mathrm{Gal}}
\def\rk{\mathrm{rank}}
\def\cork{\mathrm{corank}}

\def\div{\mathrm{div}}

\def\unr{\mathrm{ur}}

\def\Res{\mathrm{Res}}

\def\End{\mathrm{End}}

\def\Tr{\mathrm{Tr}}

\def\too{\longrightarrow}
\def\Sel{\mathrm{Sel}}
\def\Scp{\cS_p}
\def\map#1{\;\xrightarrow{#1}\;}
\def\isom{\xrightarrow{\sim}}
\def\hookto{\hookrightarrow}
\def\onto{\twoheadrightarrow}

\def\c{\mathbf{c}}

\def\Hs#1{H^1_{#1}}
\def\HF{\Hs{\cF}}
\def\HE{\Hs{\E}}
\def\HA{\Hs{\A}}
\def\HG{\Hs{\cG}}

\def\bmu{\boldsymbol{\mu}}
\def\Ehat{\hat{E}}

\def\Qpp{\Q_{p^2}}
\def\Zpp{\Z_{p^2}}
\def\dirsum#1{\underset{#1}{\textstyle\bigoplus}}
\def\alg{X}
\def\Shmd{\Sh_{/\div}}

\def\k{k}
\def\res{\kappa}

\def\rf{\Fp}

\def\Rp{R_p}
\def\Dp{D_p}

\def\Kantip{\K^-_{c,p}}
\def\semidirect{\rtimes}
\def\Sbad{\mathfrak{S}}

\def\pair#1#2{\ld#1,#2\rd}
\def\spair#1#2{[#1,#2]}
\def\bpair#1#2{\{#1,#2\}}

\title[Finding large Selmer rank]
{Finding large Selmer rank via an arithmetic theory of local constants}

\author{Barry Mazur}
\address{Department of Mathematics,
Harvard University,
Cambridge, MA 02138 USA}
\email{mazur\char`\@math.harvard.edu}

\author{Karl Rubin}
\address{Department of Mathematics,
University of California Irvine,
Irvine, CA 92697 USA}
\email{krubin\char`\@math.uci.edu}
\thanks{The authors are supported by NSF grants DMS-0403374 and DMS-0457481, respectively.}

\begin{document}

\begin{abstract}
We obtain lower bounds for Selmer ranks of elliptic curves over dihedral
extensions of number fields.

Suppose $K/\k$ is a quadratic extension of number fields, 
$E$ is an elliptic curve defined over $\k$, and $p$ is an odd prime.  
Let $\K^-$ denote the maximal abelian $p$-extension of $K$ that is unramified 
at all primes where $E$ has bad reduction and that is Galois over $\k$ 
with dihedral Galois group (i.e., the generator $c$ of 
$\Gal(K/\k)$ acts on $\Gal(\K^-/K)$ by inversion).
We prove (under mild hypotheses on $p$) that if the $\Zp$-rank of 
the pro-$p$ Selmer group $\Scp(E/K)$ is odd, then 
$\rk_{\Zp} \Scp(E/F) \ge [F:K]$ for every finite extension $F$ of $K$ in $\K^-$.
\end{abstract}

\maketitle

\setcounter{tocdepth}{1}

\section*{Introduction}
\setcounter{section}{0}
Let $K/\k$ be a quadratic extension of number fields, let $c$ be the 
nontrivial automorphism of $K/\k$, and let $E$ be an elliptic curve 
defined over $\k$.  Let $F/K$ be an abelian extension 
such that $F$ is Galois over $\k$ with 
dihedral Galois group (i.e., a lift of the involution $c$ 
operates by conjugation on $\Gal(F/K)$ as inversion $x \mapsto x^{-1}$), and 
let $\chi : \Gal(F/K) \to \bar{\Q}^\times$ be a character.  

Even in cases where one cannot prove 
that the $L$-function $L(E/K,\chi; s)$ has an analytic continuation 
and functional equation, one still has a conjectural functional equation 
with a sign $\epsilon(E/K,\chi) := \prod_v \epsilon(E/K_v,\chi_v) = \pm 1$ expressed as 
a product over places $v$ of $K$ of local $\epsilon$-factors.  
If $\epsilon(E/K,\chi) = -1$, then a generalized Parity Conjecture predicts that the 
rank of the $\chi$-part $E(F)^\chi$ 
of the $\Gal(F/K)$-representation space $E(F)\otimes \bar{\Q}$ 
is odd, and hence positive.  
If $[F:K]$ is odd and 
$F/K$ is unramified at all primes where $E$ has bad reduction, 
then $\epsilon(E/K,\chi)$ is independent of $\chi$, and so 
the Parity Conjecture predicts that 
if the rank of $E(K)$ is odd then the rank of $E(F)$ is at least $[F:K]$.

Motivated by the analytic theory of the preceding paragraph, 
in this paper we prove unconditional parity 
statements, not for the Mordell-Weil groups $E(F)^\chi$ but instead 
for the corresponding pro-$p$ Selmer groups $\Scp(E/F)^\chi$.  
(The Shafarevich-Tate conjecture implies that $E(F)^\chi$ and $\Scp(E/F)^\chi$ 
have the same rank).
More specifically, given the data $(E,K/\k,\chi)$ where the order of 
$\chi$ is a power of an odd prime $p$, 
we define (by cohomological methods) local invariants $\delta_v \in \Z/2\Z$ 
for the finite places $v$ of $K$, depending only on $E/K_v$ and 
$\chi_v$.  The $\delta_v$ should be the 
(additive) counterparts of the ratios 
$\epsilon(E/K_v,\chi_v)/\epsilon(E/K_v,1)$ of the local $\epsilon$-factors.  
The $\delta_v$ vanish for almost all $v$, and if 
$\Zp[\chi]$ is the extension of $\Zp$ generated by the values of $\chi$, 
we prove (see Theorem \ref{rkdif}):

\begin{introthm}
\label{thma}
If the order of $\chi$ is a power of an odd prime $p$, then 
$$
\rk_{\Zp}\Scp(E/K) - \rk_{\Zp[\chi]}\Scp(E/F)^\chi 
    \equiv \sum_v \delta_v \pmod{2}.
$$
\end{introthm}

Despite the fact that the analytic theory, which is our guide, 
predicts the values of the local terms $\delta_v$, Theorem \ref{thma} 
would be of limited use if we could not actually compute the 
$\delta_v$'s.  
We compute the $\delta_v$'s in substantial generality in \S\ref{loc} and \S\ref{diex}.  
This leads to our main result (Theorem \ref{third}), which we illustrate 
here with a weaker version.

\begin{introthm}
\label{thmb}
Suppose that $p$ is an odd prime, $[F:K]$ is a power of $p$,
$F/K$ is unramified at all primes where $E$ has bad reduction, and 
all primes above $p$ split in $K/\k$.  If $\rk_{\Zp}\Scp(E/K)$ is odd, 
then $\rk_{\Zp[\chi]}\Scp(E/F)^\chi$ is odd for every character $\chi$ of $G$, and 
in particular $\rk_{\Zp}\Scp(E/F) \ge [F:K]$.
\end{introthm}

If $K$ is an imaginary quadratic field and $F/K$ is unramified outside 
of $p$, then Theorem \ref{thmb} is a consequence of work of 
Cornut and Vatsal \cite{cornut, vatsal}.  In those cases 
the bulk of the Selmer module comes from Heegner points.  

Nekov\'a\u{r} (\cite{nekovar} Theorem 10.7.17) 
proved Theorem \ref{thmb} in the case where 
$F$ is contained in a $\Zp$-power extension of $K$, under the assumption 
that $E$ has ordinary reduction at all primes above $p$.  We gave in 
\cite{finding} an exposition of a weaker version of Nekov\'a\u{r}'s theorem, 
as a direct application of a functional equation that arose in 
\cite{organization} (which also depends heavily on Nekov\'a\u{r}'s theory 
in \cite{nekovar}). 

The proofs of Theorems \ref{thma} and \ref{thmb} proceed by methods 
that are very different from those of Cornut, Vatsal, and Nekov\'a\u{r}, 
and are comparatively short.  We emphasize that our results 
apply whether $E$ has ordinary or supersingular reduction at $p$, and they apply 
even when $F/K$ is not contained in a $\Zp$-power extension of $K$ 
(but we always assume that $F/\k$ is dihedral).  

This extra generality is of particular interest in connection with 
the search for new Euler systems, beyond the known examples 
of Heegner points.  Let $\K^- = \Kantip$ be the maximal 
``generalized dihedral'' $p$-extension of $K$ (i.e., the maximal abelian 
$p$-extension of $K$, Galois over $\k$, such that $c$ acts on $\Gal(\K^-/K)$ 
by inversion).  
A ``dihedral'' Euler system $\c$ for $(E,K/\k,p)$ 
would consist of Selmer classes $\c_F \in \Scp(E/F)$ for every finite 
extension $F$ of $K$ in $\K^-$, with certain compatibility relations 
between $\c_F$ and $\c_{F'}$ when $F \subset F'$ 
(see for example \cite{EulerSystems} \S9.4).  
A necessary condition for the existence of a nontrivial 
Euler system is that the Selmer modules $\Scp(E/F)$ are large, 
as in the conclusion of Theorem \ref{thmb}.  
It is natural to ask whether, in these large Selmer modules $\Scp(E/F)$, 
one can find elements $\c_F$ that form an Euler system.

\subsection*{Outline of the proofs}  
Suppose for simplicity that $E(K)$ has no $p$-torsion.

The group ring $\Q[\Gal(F/K)]$ splits into a sum of irreducible rational 
representations $\Q[\Gal(F/K)] = \oplus_L \rho_L$, 
summing over all cyclic extensions $L$ of $K$ in $F$,
where $\rho_L \otimes \bar{\Q}$ is the sum of all characters $\chi$ whose 
kernel is $\Gal(F/L)$.  Corresponding to this decomposition there is a 
decomposition (up to isogeny) of the restriction of scalars $\Res^F_K E$ 
into abelian varieties over $K$ 
$$
\Res^F_K E \sim \oplus_{L}A_L.
$$
This gives a decomposition of Selmer modules
$$
\Scp(E/F) \cong \Scp((\Res^F_K E)/K) \cong \oplus_{L}\Scp(A_L/K)
$$
where for every $L$, $\Scp(A_L/K) \cong (\rho_L \otimes \Qp)^{d_L}$ for some $d_L \ge 0$.  
Theorem \ref{thmb} will follow once we show that 
$d_L \equiv \rk_{\Zp}\Scp(E/K) \pmod{2}$ for every $L$.
More precisely, we will show (see \S\ref{localinv} for the ideal $\P$ of $\End_K(A_L)$, 
\S\ref{ecx} for the Selmer groups $\Sel_p$ and $\Sel_\P$, 
and Definition \ref{scpdef} for $\Scp$) that
\begin{equation}
\label{sketch}
\rk_{\Zp}\Scp(E/K) \equiv \dim_{\Fp}\Sel_p(E/K) 
    \equiv \dim_{\Fp}\Sel_\P(A_L/K) 
    \equiv d_L \pmod{2}.
\end{equation}

The key step in our proof 
is the second congruence of \eqref{sketch}.  
We will see (Proposition \ref{moreresprop}) 
that $E[p] \cong A_L[\P]$ as $G_K$-modules, 
and therefore the Selmer groups $\Sel_p(E/K)$ and $\Sel_\P(A_L/K)$ 
are both contained in $H^1(K,E[p])$.  By comparing these two 
subspaces we prove (see Theorem \ref{fg} and Corollary \ref{rkdif1}) that 
$$
\dim_{\Fp}\Sel_p(E/K) - \dim_{\Fp}\Sel_\P(A_L/K) \equiv \sum_v \delta_v \pmod{2}
$$
summing the local invariants $\delta_v$ of Definition \ref{deltadef} 
over primes $v$ of $K$.
We show how to compute the $\delta_v$ in terms of 
norm indices in \S\ref{loc} and \S\ref{diex}, with one important 
special case postponed to Appendix \ref{appb}.

The first congruence of \eqref{sketch} follows easily from the Cassels 
pairing for $E$ (see Proposition \ref{selfA}).   
The final congruence of \eqref{sketch} is more subtle, because 
in general $A_L$ will not have a polarization of degree prime to $p$, and 
we deal with this in Appendix \ref{polpf} (using the dihedral nature of $L/\k$).

In \S\ref{d1} we bring together the results of the previous sections 
to prove Theorem \ref{third}, and in 
\S\ref{examples} we discuss some special cases.   

\subsection*{Generalizations}
All the results and proofs in this paper hold with $E$ replaced by an 
abelian variety with a polarization of degree prime to $p$.  

If $F/K$ is not a $p$-extension, then the proof described above breaks down.  
Namely, if $\chi$ is a character whose order is not a prime power, 
then $\chi$ is not congruent to the trivial character modulo any prime 
of $\bar\Q$.  However, by writing $\chi$ as a product of characters 
of prime-power order, we can apply the methods of this paper 
inductively.  To do this we must use a different prime $p$ at each step, 
so it is necessary to assume that if $A$ is an abelian variety over $K$ 
and $R$ is an integral domain in 
$\End_K(A)$, then the parity of $\dim_{R \otimes \Qp}\Scp(A/K)$ is 
independent of $p$.  (This would follow, for example, from the Shafarevich-Tate 
conjecture.)  
To avoid obscuring the main ideas of our arguments, 
we will include those details in a separate paper.

The results of this paper can also be applied to study the growth of Selmer 
rank in nonabelian Galois extensions of order $2p^n$ with $p$ an odd prime.  
This will be the subject of a forthcoming paper.

\subsection*{Notation}  
Fix once and for all an algebraic closure $\bar{\Q}$ of $\Q$.  
A number field will mean a finite extension of $\Q$ in $\bar{\Q}$.  
If $K$ is a number field then $G_K := \Gal(\bar{\Q}/K)$.

\numberwithin{equation}{section}

\section{Variation of Selmer rank}
\label{var}

Let $K$ be a number field and $p$ an odd rational prime.   
Let $\T$ be a finite-dimensional $\rf$-vector space with a continuous 
action of $G_K$ and with a perfect, skew-symmetric, 
$G_K$-equivariant self-duality 
$$
\T \times \T \too  \bmu_p 
$$
where $\bmu_p$ is the $G_K$-module of $p$-th roots of unity in $\bar{\Q}$.  

\begin{thm}
\label{tld}
For every prime $v$ of $K$, 
Tate's local duality gives a perfect symmetric pairing
$$
\pair{\;}{\;}_v : H^1(K_v,\T) \times H^1(K_v,\T) \too H^2(K_v,\bmu_p) = \Fp.
$$
\end{thm}

\begin{proof}
See \cite{tate}.
\end{proof}

\begin{defn}
\label{ss}
For every prime $v$ of $K$ let 
$K_v^\unr$ denote the maximal unramified extension of $K_v$.
A {\em Selmer structure} $\cF$ on $\T$ is a collection of $\rf$-subspaces 
$$
\HF(K_v,\T) \subset H^1(K_v,\T)
$$
for every prime $v$ of $K$, such that $\HF(K_v,\T) = H^1(K_v^\unr/K_v,\T^{I_v})$ 
for all but finitely many $v$, where $I_v := G_{K_v^\unr} \subset G_{K_v}$ is the inertia group.  
If $\cF$ and $\cG$ are Selmer structures on $\T$, we define 
Selmer structures $\cF+\cG$ and $\cF\cap\cG$ by 
\begin{align*}
\Hs{\cF+\cG}(K_v,\T) &:= \HF(K_v,\T)+\HG(K_v,\T), \\
\Hs{\cF\cap\cG}(K_v,\T) &:= \HF(K_v,\T)\cap\HG(K_v,\T)
\end{align*}
for every $v$.  We say that $\cF \le \cG$ if $\HF(K_v,\T) \subset \HG(K_v,\T)$ 
for every $v$, so in particular $\cF \cap \cG \le \cF \le \cF+\cG$.

We say that a Selmer structure $\cF$ is {\em self-dual} if 
for every $v$, $\HF(K_v,\T)$ is its own orthogonal complement 
under the Tate pairing of Theorem \ref{tld}.

If $\cF$ is a Selmer structure on $\T$, we define the {\em Selmer group} 
$$
\HF(K,\T) := \ker(H^1(K,\T) \too \textstyle\prod_{v} H^1(K_v,\T)/\HF(K_v,\T)). 
$$
Thus $\HF(K,\T)$ is 
the collection of classes whose localizations lie in $\HF(K_v,\T)$ 
for every $v$.  If $\cF \le \cG$ then $\HF(K,\T) \subset \HG(K,\T)$.
\end{defn}

For the basic example of the Selmer groups we 
will be interested in, where $\T$ is the Galois module of 
$p$-torsion on an elliptic curve, see \S\ref{ecx}.

\begin{prop}
\label{globdual}
Suppose that
$\cF$, $\cG$ are self-dual Selmer structures on $\T$, 
and $S$ is a finite set of primes of $K$ 
such that $\HF(K_v,\T) = \HG(K_v,\T)$ if $v \notin S$.  Then
\begin{enumerate}
\item
$\dim_{\rf}\Hs{\cF+\cG}(K,\T)/\Hs{\cF\cap\cG}(K,\T) 
    = \displaystyle\sum_{v \in S}\dim_{\rf}\HF(K_v,\T)/\Hs{\cF\cap\cG}(K_v,\T)$,
\item
$\dim_{\rf}\Hs{\cF+\cG}(K,\T) \equiv \dim_{\rf}(\HF(K,\T)+\HG(K,\T)) \pmod{2}$.
\end{enumerate}
\end{prop}

\begin{proof}
Let 
$$
B := \dirsum{v \in S}(\Hs{\cF+\cG}(K_v,\T)/\Hs{\cF\cap\cG}(K_v,\T))
$$ 
and let $C$ be the image of the localization map $\Hs{\cF+\cG}(K,\T) \to B$.
Since $\cF$ and $\cG$ are self-dual, 
Poitou-Tate global duality (see for example \cite{kolysys} Theorem 2.3.4)
shows that 
the Tate pairings of Theorem \ref{tld} induce a nondegenerate, symmetric self-pairing 
\begin{equation}
\label{derlp}
\pair{\;\,}{\;} : B \times B \too \Fp,
\end{equation}
and $C$ is its own orthogonal complement under this pairing.  

Let $C_\cF$ (resp.\ $C_\cG$) denote the image of $\oplus_{v \in S}\HF(K_v,\T)$ 
(resp.\ $\oplus_{v \in S}\HG(K_v,\T)$) in $B$.
Since $\cF$ and $\cG$ are self-dual, $C_\cF$ and $C_\cG$ are 
each their own orthogonal complements under \eqref{derlp}.  
In particular we have
$$
\dim_{\rf}C = \dim_{\rf}C_\cF = \dim_{\rf}C_\cG = \textstyle\frac{1}{2}\dim_{\rf}B.
$$
Since $C \cong \Hs{\cF+\cG}(K,\T)/\Hs{\cF\cap\cG}(K,\T)$ and 
$C_\cF \cong \oplus_{v \in S}\HF(K_v,\T)/\Hs{\cF\cap\cG}(K_v,\T)$, this proves (i).

The proof of (ii) uses an argument of Howard (\cite{howard} Lemma 1.5.7).  
We have $C_\cF \cap C_\cG = 0$ and $C_\cF \oplus C_\cG = B$.  
If $x \in \Hs{\cF+\cG}(K,\T)$, let $x_S \in C \subset B$ be the 
localization of $x$, and let $x_\cF$ and $x_\cG$ denote 
the projections of $x_S$ to $C_\cF$ and $C_\cG$, respectively.  

Following Howard, we define a pairing 
\begin{equation}
\label{rederlp}
[\;\,,\;] : \Hs{\cF+\cG}(K,\T) \times \Hs{\cF+\cG}(K,\T) \too \Fp
\end{equation}
by $[x,y] := \pair{x_\cF}{y_\cG}$, where $\pair{\;\,}{\;}$ is the pairing 
\eqref{derlp}.
Since the subspaces $C$, $C_\cF$, and $C_\cG$ are all isotropic, 
for all $x, y, \in \Hs{\cF+\cG}(K,\T)$ we have
$$
0 = \pair{x_S}{y_S} = \pair{x_\cF+x_\cG}{y_\cF+y_\cG} 
    = \pair{x_\cF}{y_\cG} + \pair{x_\cG}{y_\cF} 
    = [x,y] + [y,x]
$$
so the pairing \eqref{rederlp} is skew-symmetric.  

We see easily that $\HF(K,\T) + \HG(K,\T)$ is in the kernel of the pairing $[\;\,,\;]$.  
Conversely, 
if $x$ is in the kernel of this pairing, then 
for every $y \in \Hs{\cF+\cG}(K,\T)$
$$
0 = [x,y] = \pair{x_\cF}{y_\cG} = \pair{x_\cF}{y_S}.
$$
Since $C$ is its own orthogonal complement we deduce that $x_\cF \in C$, 
i.e., there is a $z \in \Hs{\cF+\cG}(K,\T)$ 
whose localization is $x_\cF$.  It follows that $z \in \HF(K,\T)$ 
and $x-z \in \HG(K,\T)$, i.e., $x \in \HF(K,\T) + \HG(K,\T)$.  
Therefore \eqref{rederlp} induces a nondegenerate, skew-symmetric, $\Fp$-valued 
pairing on $$\Hs{\cF+\cG}(K,\T)/(\HF(K,\T) + \HG(K,\T)).$$ 
Since $p$ is odd, a well-known argument from linear algebra shows that 
the dimension of this $\Fp$-vector space must be even.  This proves (ii).
\end{proof}

\begin{thm}
\label{fg}
Suppose that $\cF$ and $\cG$ are self-dual Selmer structures on $\T$, and 
$S$ is a finite set of primes of $K$ 
such that $\HF(K_v,\T) = \HG(K_v,\T)$ if $v \notin S$.  Then
\begin{multline*}
\dim_{\rf}\HF(K,\T) - \dim_{\rf}\HG(K,\T) \\ \equiv \sum_{v \in S} 
    \dim_{\rf}(\HF(K_v,\T)/\Hs{\cF\cap\cG}(K_v,\T)) \pmod{2}.
\end{multline*}
\end{thm}

\begin{proof}
We have (modulo $2$)
\begin{align*}
\dim_{\rf} \HF(K,\T) - &\dim_{\rf} \HG(K,\T) 
    \equiv \dim_{\rf} \HF(K,\T) + \dim_{\rf} \HG(K,\T) \\
    &= \dim_{\rf} (\HF(K,\T) + \HG(K,\T)) + \dim_{\rf} \Hs{\cF\cap\cG}(K,\T) \\
    &\equiv \dim_{\rf}\Hs{\cF+\cG}(K,\T) - \dim_{\rf}\Hs{\cF\cap\cG}(K,\T) \\
    &=  \sum_{v \in S}\dim_{\rf}(\HF(K_v,\T)/\Hs{\cF\cap\cG}(K_v,\T)),
\end{align*}
the last two steps by Proposition \ref{globdual}(ii) and (i), respectively.
\end{proof}

\section{Example: elliptic curves}
\label{ecx}
Let $K$ be a number field.
If $A$ is an abelian variety over $K$, and $\alpha \in \End_K(A)$ is 
an isogeny, we have the usual Selmer group $\Sel_\alpha(A/K) \subset H^1(K,E[\alpha])$, 
sitting in an exact sequence  
\begin{equation}
\label{selE}
0 \too A(K)/\alpha A(K) \too \Sel_\alpha(A/K) \too \Sh(A/K)[\alpha] \too 0
\end{equation}
where $\Sh(A/K)$ is the Shafarevich-Tate group of $A$ over $K$.
If $p$ is a prime we let $\Sel_{p^\infty}(A/K)$ be the direct limit of the 
Selmer groups $\Sel_{p^n}(A/K)$, and then we have
\begin{equation}
\label{selEinf}
0 \too A(K) \otimes \Qp/\Zp \too \Sel_{p^\infty}(A/K) 
    \too \Sh(A/K)[p^\infty] \too 0.
\end{equation}

Suppose now that $E$ is an elliptic curve defined over $K$, 
and $p$ is an odd rational prime.  
Let $\T := E[p]$, the Galois module of $p$-torsion in $E(\bar{\Q})$.  
Then $\T$ is an $\rf$-vector space with a continuous action of $G_K$, 
and the Weil pairing induces a perfect $G_K$-equivariant self-duality 
$E[p] \times E[p] \to \bmu_p$.  Thus we are in the setting of \S\ref{var}.

We define a Selmer structure $\E$ on $E[p]$ by taking 
$\HE(K_v,E[p])$ to be the image of $E(K_v)/p E(K_v)$ under the Kummer injection
$$
E(K_v)/p E(K_v) \hookto H^1(K_v, E[p])
$$
for every $v$.  By Lemma 19.3 of \cite{casexp}, 
$\HE(K_v,E[p]) = H^1(K_v^\unr/K_v,E[p])$ if $v \nmid p$ and 
$E$ has good reduction at $v$.
With this definition the Selmer group $\HE(K,E[p])$ is 
the usual $p$-Selmer group $\Sel_p(E/K)$ of $E$ as in \eqref{selE}.

If $C$ is an abelian group, we let $C_\div$ denote its maximal divisible subgroup.  

\begin{prop}
\label{selfA}
The Selmer structure $\E$ on $E[p]$ defined above is self-dual, and
$$
\cork_{\Zp}\Sel_{p^\infty}(E/K) \equiv \dim_{\rf}\HE(K,E[p]) - \dim_{\rf}E(K)[p] \pmod{2}.
$$  
\end{prop}

\begin{proof}
Tate's local duality \cite{tate} shows that $\E$ is self-dual.  
Let 
\begin{align*}
d &:= \dim_{\rf}(\Sel_{p^\infty}(E/K)/(\Sel_{p^\infty}(E/K))_\div)[p] \\
    &\phantom{:}= \dim_{\rf}(\Sh(E/K)[p^\infty]/(\Sh(E/K)[p^\infty])_\div)[p].
\end{align*}
The Cassels pairing \cite{cassels} shows that $d$ is even.  
Further, 
\begin{align*}
\cork_{\Zp}\Sel_{p^\infty}(E/K) &= \dim_{\rf}\Sel_{p^\infty}(E/K)_\div[p] \\
    &= \dim_{\rf}\Sel_{p^\infty}(E/K)[p] - d \\
    &= \rk_\Z E(K) + \dim_{\rf}\Sh(E/K)[p] - d
\end{align*}
by \eqref{selEinf} with $A = E$.  On the other hand, \eqref{selE} shows that 
$$
\dim_{\rf}\HE(K,E[p]) = \rk_\Z E(K) + \dim_{\rf}E(K)[p] + \dim_{\rf}\Sh(E/K)[p]
$$
so we conclude
$$
\cork_{\Zp}\Sel_{p^\infty}(E/K) = \dim_{\rf}\HE(K,E[p]) - \dim_{\rf}E(K)[p] - d.
$$
This proves the proposition.
\end{proof}

\section{Decomposition of the restriction of scalars}
\label{resn}

Much of the technical machinery for this section will be drawn 
from sections 4 and 5 of \cite{prim}.

Suppose $F/K$ is a finite abelian extension of number fields, 
$G := \Gal(F/K)$, and $E$ is an elliptic curve defined over $K$.  
We let $\Res^F_K E$ denote the Weil restriction of scalars (\cite{weil} \S1.3)
of $E$ from $F$ to $K$, an abelian variety over $K$ 
with the following properties.

\begin{prop}
\label{newresprop}
\begin{enumerate}
\item
For every commutative $K$-algebra $\alg$ there is a 
canonical isomorphism 
$$
(\Res^F_K E)(\alg) \cong E(\alg \otimes_K F)
$$
functorial in $\alg$.  In particular $(\Res^F_K E)(K) \cong E(F)$.
\item
The action of $G$ on the right-hand side of (i) induces a canonical 
inclusion $\Z[G] \hookto \End_K(\Res^F_K E)$.
\item
For every prime $p$ there is a natural $G$-equivariant isomorphism, 
compatible with the isomorphism $(\Res^F_K E)(K) \cong E(F)$ of (i),
$$
\Sel_{p^\infty}((\Res^F_K E)/K) \cong \Sel_{p^\infty}(E/F)
$$
where $G$ acts on the left-hand side via the inclusion of (ii). 
\end{enumerate}
\end{prop}

\begin{proof}
Assertion (i) is the universal property satisfied by the restriction 
of scalars \cite{weil}, and (ii) is (for example) (4.2) of \cite{prim}.  
For (iii), Theorem 2.2(ii) and Proposition 4.1 of \cite{prim} give an isomorphism  
$$
(\Res^F_K E)[p^\infty] \cong \Z[G] \otimes E[p^\infty]
$$ 
that is $G$-equivariant (with $G$ acting on $\Res^F_K E$ via the map of (ii) 
and by multiplication on $\Z[G]$) and $G_K$-equivariant 
(with $\gamma \in G_K$ acting by 
$\gamma^{-1} \otimes \gamma$ on $\Z[G] \otimes E[p^\infty]$).
Then by Shapiro's Lemma 
(see for example Propositions III.6.2, III.5.6(a), and III.5.9 of \cite{brown}) 
there is a $G$-equivariant isomorphism
\begin{equation}
\label{resiso}
H^1(K,(\Res^F_K E)[p^\infty]) \isom H^1(F,E[p^\infty]).
\end{equation}
Using (i) with $\alg = K_v$, along with the analogue of \eqref{resiso}
for the local extensions $(F \otimes_K K_v)/K_v$ for every 
prime $v$ of $K$, one can show that the isomorphism \eqref{resiso} restricts
to the isomorphism of (iii).
\end{proof}

\begin{defn}
\label{xirndef}
Let $\Xi := \{\text{cyclic extensions of $K$ in $F$}\}$, and if $L \in \Xi$ 
let $\rho_L$ be the unique faithful irreducible rational representation of 
$\Gal(L/K)$.  Then $\rho_L \otimes \bar\Q$ is the 
direct sum of all the injective characters $\Gal(L/K) \hookto \bar{\Q}^\times$.  
The correspondence $L \leftrightarrow \rho_L$ is a bijection between 
$\Xi$ and the set of irreducible rational representations of $G$.
Thus the semisimple group ring $\Q[G]$ decomposes 
\begin{equation}
\label{ssgd}
\Q[G] \cong \dirsum{L \in \Xi} \Q[G]_L
\end{equation}
where $\Q[G]_L \cong \rho_L$ is the $\rho_L$-isotypic component of $\Q[G]$.  
As a field, $\Q[G]_L$ is isomorphic to the cyclotomic field of $[L:K]$-th 
roots of unity.  

Let $R_L$ be the maximal order of $\Q[G]_L$.  If $[L:K]$ is a power 
of a prime $p$, then $R_L$ has a unique prime ideal above $p$, which we 
denote by $\P_L$.  
Also define
$$
\I_L := \Q[G]_L \cap \Z[G], 
$$
so $\I_L$ is an ideal of $R_L$ as well as a $G_K$-module, 
(where the action of $G_K$ is induced by multiplication on $\Z[G]$).
\end{defn}

\begin{defn}
\label{ildef}
For every $L \in \Xi$ define 
$$
A_L := \I_L \otimes E
$$
as given by Definition 1.1 of \cite{prim}  
(see also \cite{milne} \S2).  
The abelian variety $A_L$ is defined over $K$, and 
its $K$-isomorphism class is independent of the choice of 
abelian extension $F$ containing $L$ 
(see Remark 4.4 of \cite{prim}).   
If $L = K$ then $A_K = E$.
By Proposition 4.2(i) of \cite{prim},  
the inclusion $\I_L \hookto \Z[G]$ induces an isomorphism 
\begin{equation}
\label{alrese}
A_L \cong \bigcap_{\alpha \in \Z[G] \;:\; \alpha \I_L = 0} 
    \ker(\alpha : \Res^F_K E \to \Res^F_K E) ~\subset~ \Res^F_K E.
\end{equation}
\end{defn}

Let $T_p(E)$ denote the Tate module $\varprojlim E[p^n]$, and similarly 
for $T_p(A_L)$.
The following theorem summarizes the properties of the abelian varieties $A_L$ 
that we will need.

\begin{thm}
\label{athm}
Suppose $p$ is a prime, $n \ge 1$, and 
$L/K$ is a cyclic extension of degree $p^n$.  Then:
\begin{enumerate}
\item
$\I_L = \P_L^{p^{n-1}}$ in $R_L$.
\item
The inclusion 
$\Z[G] \hookto \End_K(\Res^F_K E)$ of Proposition \ref{newresprop}(ii)
induces (via \eqref{alrese}) a ring homomorphism $\Z[G] \to \End_K(A_L)$ 
that factors 
$$
\Z[G] \onto R_L \hookto \End_K(A_L)
$$
where the first map is induced by the projection in \eqref{ssgd}.
\item
Let $M$ be the unique extension of $K$ in $L$ with $[L:M] = p$.  
For every commutative $K$-algebra $\alg$, the isomorphism of 
Proposition \ref{newresprop}(i) restricts (using \eqref{alrese}) to an isomorphism, 
functorial in $\alg$, 
$$
A_L(\alg) \cong \{x \in E(\alg \otimes_K L) : 
    \sum_{h \in \Gal(L/M)}(1 \otimes h)(x) = 0\}.
$$
\item
The isomorphism of (iii) with $X = \bar{\Q}$ induces an isomorphism
$$
T_p(A_L) \cong \I_L \otimes T_p(E) = \P_L^{p^{n-1}} \otimes T_p(E)
$$
that is $G_K$-equivariant, where $\gamma \in G_K$ acts on 
the tensor products as $\gamma^{-1} \otimes \gamma$, 
and $R_L$-linear, where $R_L$ acts on $A_L$ via the map of (ii).
\end{enumerate}
\end{thm}

\begin{proof}
Assertions (i), (ii), and (iv) are Lemma 5.4(iv), Theorem 5.5(iv), 
and Theorem 2.2(iii), respectively, of \cite{prim} 
((iv) is also Proposition 6(b) of \cite{milne}).  
Assertion (iii) is Theorem 5.8(ii) of \cite{prim}.  
\end{proof}

\begin{thm}
\label{newdecomp}
The inclusions $A_L \subset \Res^F_K E$ of \eqref{alrese} induce an isogeny 
$$
\dirsum{L \in \Xi} A_L \too \Res^F_K E.
$$
\end{thm}

\begin{proof}
This is Theorem 5.2 of \cite{prim};  
it follows from the fact that $\oplus_{L \in \Xi} \;\I_L$ injects into 
$\Z[G]$ with finite cokernel.
\end{proof}

\begin{defn}
\label{scpdef}
Define the Pontrjagin dual Selmer vector spaces 
\begin{align*}
\Scp(E/K) &:= \Hom(\Sel_{p^\infty}(E/K),\Qp/\Zp) \otimes \Qp, \\
\Scp(A_L/K) &:= \Hom(\Sel_{p^\infty}(A_L/K),\Qp/\Zp) \otimes \Qp.
\end{align*}
Define $\Scp(E/F)$ similarly for every finite extension $F$ of $K$.
\end{defn}

\begin{cor}
\label{newressca}
There is a $G$-equivariant isomorphism 
$$
\Scp(E/F) \cong \dirsum{L\in\Xi} \Scp(A_L/K)
$$
where the action of $G$ on the right-hand side is 
given by Theorem \ref{athm}(ii).
\end{cor}

\begin{proof}
We have $\Scp(E/F) \cong \Scp((\Res^F_K E)/K)$ by (the Pontrjagin dual of) 
Proposition \ref{newresprop}(iii), and 
$\Scp((\Res^F_K E)/K) \cong \oplus_{L \in \Xi} \Scp(A_L/K)$
by Theorem \ref{newdecomp}.
\end{proof}

\section{The local invariants}
\label{localinv}

Fix an odd prime $p$ and a cyclic extension $L/K$ of degree $p^n$.  We will 
write simply $A$ for the abelian variety $A_L$ of Definition \ref{ildef},  
$R$ for the ring $R_L$ of Definition \ref{xirndef}, $\P$ for the unique prime 
$\P_L$ of $R$ above $p$, and $\I \subset R$ for the ideal $\I_L$ of 
Definition \ref{xirndef}.

\begin{prop}
\label{moreresprop}
There is a canonical $G_K$-isomorphism $A[\P] \isom E[p]$.
\end{prop}

\begin{proof}
The action of $G$ on $\P^{-1}\I/\I$ is trivial, 
since for every $g \in G$, $g-1$ lies in the maximal ideal of $\Zp[G]$.
Also, if $\pi$ and $\pi'$ are generators 
of $\P/\P^2$, then $\pi/\pi' \in (R/\P)^\times = \Fp^\times$, so 
$\pi^{p-1} \equiv (\pi')^{p-1} \pmod{\P^{p}}$.  It follows that 
$\pi^{p-1}$ is a canonical generator of $\P^{p-1}/\P^p$, so there is 
a canonical isomorphism $\P^{a(p-1)}/\P^{a(p-1)+1} \cong \Fp$ 
for every integer $a$.
Now using Theorem \ref{athm}(iv) we have $G_K$-isomorphisms 
$$
A[\P] \cong \P^{-1}T_p(A)/T_p(A) \cong (\P^{p^{n-1}-1}/\P^{p^{n-1}}) \otimes T_p(E) 
    \cong\Fp \otimes T_p(E) \cong E[p].
$$
\end{proof}

\begin{rem}
Identifying $E$ with $A_K$, one can show using \eqref{alrese} 
that $$E[p] = E \cap A_L = A_L[\P]$$ inside $\Res^F_K E$.  
This gives an alternate proof of Proposition \ref{moreresprop}.
\end{rem}

\begin{defn}
\label{bdef}
Recall that in \S\ref{ecx} we defined a self-dual Selmer structure $\E$ on $E[p]$.  
We can use the identification of Proposition \ref{moreresprop} 
to define another Selmer structure $\A$ on $E[p]$ as follows.  
For every $v$ define 
$\HA(K_v,E[p])$ to be the image of $A(K_v)/\P A(K_v)$ under the 
composition
\begin{equation*}
\label{locA}
A(K_v)/\P A(K_v) \hookto H^1(K_v,A[\P]) \cong H^1(K_v,E[p])
\end{equation*}
where the first map is the Kummer injection, 
and the second map is from Proposition \ref{moreresprop}.  
The first map depends (only up to multiplication by a unit in $\Fp^\times$) 
on a choice of generator of $\P/\P^2$, 
but the image is independent of this choice.  
With this definition the Selmer group $\HA(K,E[p])$ is 
the usual $\P$-Selmer group $\Sel_\P(A/K)$ of $A$, as in \eqref{selE}.
\end{defn}

\begin{prop}
\label{selfB}
The Selmer structure $\A$ is self-dual.
\end{prop}

\begin{proof}
This is Proposition \ref{sd} of Appendix \ref{polpf}.  
(It does not follow immediately from Tate's local duality 
as in Proposition \ref{selfA}, because $A$ has no polarization 
of degree prime to $p$, and hence no suitable Weil pairing.)
\end{proof}

\begin{defn}
\label{deltadef}
For every prime $v$ of $K$ we define an invariant $\delta_v \in \Z/2\Z$ by 
$$
\delta_v = \delta(v,E,L/K) 
    := \dim_{\rf}(\HE(K_v,E[p])/\Hs{\E\cap\A}(K_v,E[p])) \pmod{2}.
$$
We will see in Corollary \ref{intcor} below that $\delta_v$ is a purely 
local invariant, depending only on $K_v$, $E/K_v$, and $L_w$, where $w$ is a 
prime of $L$ above $v$.
\end{defn}

\begin{cor}
\label{rkdif1}
Suppose that $S$ is a set of primes of $K$ containing all primes above $p$, 
all primes ramified in $L/K$, and all primes where $E$ has bad 
reduction.  Then
$$
\dim_{\rf}\Sel_p(E/K) - \dim_{\rf}\Sel_\P(A/K) 
   \equiv \sum_{v \in S} \delta_v \pmod{2}.
$$
\end{cor}

\begin{proof}
If $v \notin S$ then both $T_p(E)$ and $T_p(A)$ are unramified at $v$, 
so (see for example \cite{casexp} Lemma 19.3) 
$$
\HE(K_v,E[p]) = \HA(K_v,E[p]) = H^1(K_v^\unr/K_v,E[p]). 
$$
Thus the corollary follows from Propositions \ref{selfA} and \ref{selfB}  
and Theorem \ref{fg}.
\end{proof}

\section{Computing the local invariants}
\label{loc}

Let $p$, $L/K$, $A := A_L$, and $\p \subset R$ be as in \S\ref{localinv}.  
Let $M$ be the unique extension of $K$ in $L$ with $[L:M] = p$, and
let $G := \Gal(L/K)$ (recall that $L/K$ is cyclic of degree $p^n$).  
In this section we compare the local Selmer conditions 
$\HE(K_v,E[p])$ and $\HA(K_v,E[p])$ for primes $v$ of $K$, in order to 
compute the invariants $\delta_v$ of Definition \ref{deltadef}.  

\begin{lem}
\label{par2}
Suppose that $c$ is an automorphism of $K$, and $E$ is defined over 
the fixed field of $c$ in $K$.  
Then for every prime $v$ of $K$, we have $\delta_{v^c} = \delta_v$.
\end{lem}

\begin{proof}
The automorphism $c$ induces isomorphisms 
$$
E(K_v) \isom E(K_{v^c}), \quad A(K_v) \isom A(K_{v^c}).
$$ 
Therefore the 
isomorphism $H^1(K_v,E[p]) \isom H^1(K_{v^c},E[p])$ induced by $c$ identifies 
$$
\HE(K_v,E[p]) \isom \HE(K_{v^c},E[p]), \quad \HA(K_v,E[p]) \isom \HA(K_{v^c},E[p]),
$$
and the lemma follows directly from the definition of $\delta_v$.
\end{proof}

For every prime $v$ of $K$, 
let $L_v := K_v \otimes_K L = \oplus_{w \mid v} L_w$, and let $G:=\Gal(L/K)$ act on 
$L_v$ via its action on $L$.  Let $M_v := K_v \otimes M$ and let 
$N_{L/M} : E(L_v) \to E(M_v)$ denote the norm (or trace) map.  
The following is our main tool for computing $\delta_v$.

\begin{prop}
\label{newint}
For every prime $v$ of $K$, the isomorphism $$\HE(K_v,E[p]) \cong E(K_v)/p E(K_v)$$ 
identifies 
$$
\Hs{\E\cap\A}(K_v,E[p]) \cong (E(K_v) \cap N_{L/M}E(L_v))/p E(K_v).
$$
\end{prop}

\begin{proof}
Fix a generator $\sigma$ of $G$, and let $\pi$ be the projection of $\sigma-1$ 
to $R$ under \eqref{ssgd}.  Since $\sigma$ projects to a $p^n$-th root of 
unity in $R$, we see that $\pi$ is a generator of $\P$.

Note that $G$ and $G_{K_v}$ act on $E(\bar{K_v} \otimes L)$ (as 
$1 \otimes G$ and $G_{K_v} \otimes 1$, respectively).  
We identify $E(L_v)$, $E(\bar{K_v})$, $A(K_v)$, and $A(\bar{K_v})$ 
with their images in $E(\bar{K_v} \otimes L)$ under the natural inclusions 
and Theorem \ref{athm}(iii):
\begin{gather*}
A(K_v) \subset E(L_v) = E(K_v \otimes L) = E(\bar{K_v} \otimes L)^{G_{K_v}}, \\
E(\bar{K_v}) = E(\bar{K_v} \otimes K) = E(\bar{K_v} \otimes L)^G, 
    \quad A(\bar{K_v}) \subset E(\bar{K_v} \otimes L).
\end{gather*}
Let $\hat\pi := (1 \otimes\sigma)-1$ on $E(\bar{K_v} \otimes L)$, 
so $\hat\pi$ restricts to $\pi$ on $A(\bar{K_v})$ and to zero on $E(\bar{K_v})$.
By Proposition \ref{athm}(iii), $A(\bar{K_v})$ is the kernel of 
$N_{L/M} := \sum_{g \in \Gal(L/M)} 1 \otimes g$ in $E(\bar{K_v} \otimes L)$.

If $x \in E(K_v)$, then the image of $x$ in $H^1(K_v,E[p])$ is represented 
by the cocycle $\gamma \mapsto y^{\gamma \otimes 1} - y$ where $y \in E(\bar{K_v})$ 
and $py = x$.  Similarly, using the identifications above, 
if $\alpha \in A(K_v)$ then the image of 
$\alpha$ in $H^1(K_v,E[p])$ is represented 
by the cocycle $\gamma \mapsto \beta^{\gamma \otimes 1} - \beta$ where 
$\beta \in A(\bar{K_v})$ and $\pi\beta = \alpha$.

Suppose $x \in E(K_v)$, and choose $y \in E(\bar{K_v})$ such that $py = x$.  Then 
\begin{align*}
\text{the i}&\text{mage of $x$ in $\HE(K_v,E[p]) \subset H^1(K_v,E[p])$ 
    belongs to $\Hs{\E\cap\A}(K_v,E[p])$} \\
&\iff \exists \beta \in A(\bar{K_v}): \pi\beta \in A(K_v), 
    \beta^{\gamma \otimes 1} - \beta = y^{\gamma \otimes 1} - y  \;\;\;\forall \gamma \in G_{K_v}\\
&\iff \exists \beta \in A(\bar{K_v}):  
    \beta^{\gamma \otimes 1} - \beta = y^{\gamma \otimes 1} - y  \;\;\;\forall \gamma \in G_{K_v}\\
&\iff \exists \beta \in E(\bar{K_v} \otimes L): N_{L/M}\beta=0,
    y-\beta \in E(L_v) \\
&\iff N_{L/M} y \in N_{L/M}E(L_v)
\end{align*}
where for the second equivalence we use that if $\gamma \in G_{K_v}$ 
and $\beta^{\gamma \otimes 1} - \beta = y^{\gamma \otimes 1} - y$, 
then $\hat\pi\beta^{\gamma \otimes 1} - \hat\pi\beta = \hat\pi(y^{\gamma \otimes 1}-y) = 0$, 
and if this holds for every $\gamma$ then $\pi\beta \in A(K_v)$.
Since $y \in E(\bar{K_v}) = E(\bar{K_v} \otimes L)^G$, we have 
$N_{L/M}y = py = x$ and the proposition follows.
\end{proof}

The following corollary gives a purely local formula for $\delta_v$, depending only on 
$E$ and the local extension $L_w/K_v$ (where $w$ is a prime of $L$ above $v$).

\begin{cor}
\label{intcor}
Suppose $v$ is a prime of $K$ and $w$ is a prime of $L$ above $v$.  If $L_w \ne K_v$ 
then let $L_w'$ be the unique subfield of $L_w$ containing $K_v$ with 
$[L_w:L_w'] = p$, and otherwise let $L_w' := L_w = K_v$.  Let 
$N_{L_w/L_w'}$ denote the norm map $E(L_w) \to E(L_w')$.  
Then
$$
\delta_v \equiv \dim_{\rf}E(K_v)/(E(K_v) \cap N_{L_w/L_w'}E(L_w)) \pmod{2}.
$$
In particular if $N_{L_w/L_w'} : E(L_w) \to E(L_w')$ is surjective 
(for example, if $v$ splits completely in $L/K$) then $\delta_v = 0$.
\end{cor}

\begin{proof}
By Proposition \ref{newint} 
$$
\HE(K_v,E[p])/\Hs{\E\cap\A}(K_v,E[p]) \cong E(K_v)/(E(K_v) \cap N_{L/M}E(L_v)), 
$$
and $\delta_v$ is the $\rf$-dimension (modulo $2$) of the left-hand side.  
Since $L/K$ is cyclic, 
$L_w'$ is the completion of $M$ at the prime below $w$, so we have
$$
E(K_v) \cap N_{L/M}E(L_v) = E(K_v) \cap N_{L_w/L_w'}E(L_w).
$$
This proves the corollary.
\end{proof}

By {\em local field} we mean a finite extension of $\Q_\ell$ for some 
rational prime $\ell$.

\begin{lem}
\label{parities}
If $\K$ is a local field with residue characteristic different from $p$, 
and $E$ is defined over $\K$, 
then $E(\K)/pE(\K) = E(\K)[p^\infty]/p E(\K)[p^\infty]$ and in particular
$$
\dim_{\rf} E(\K)/p E(\K) = \dim_{\rf}E(\K)[p].
$$
\end{lem}

\begin{proof}
There is an isomorphism of topological groups
$$
E(\K) \cong E(\K)[p^\infty] \oplus C \oplus D
$$
with a finite group $C$ of order prime to $p$ and 
a free $\Z_\ell$ module $D$ of finite rank, where $\ell$ is the residue 
characteristic of $v$.  Since $E(\K)[p^\infty]$ is finite, the lemma follows easily.  
\end{proof}

\begin{lem}
\label{newnormlem}
Suppose $\L/\K$ is a cyclic extension of degree $p$ of local fields 
and $E$ is defined over $\K$.  
Let $\ell$ denote the residue characteristic of $\K$.
\begin{enumerate}
\item
If $\L/\K$ is unramified and $E$ has good reduction, 
then $N_{\L/\K}E(\L) = E(\K)$.
\item
If $\L/\K$ is ramified, $\ell \ne p$, and 
$E$ has good reduction, then 
$$
E(\K)/p E(\K) \to E(\L)/p E(\L)
$$ 
is an isomorphism and $N_{\L/\K}E(\L) = p E(\K)$.
\end{enumerate}
\end{lem}

\begin{proof}
The first assertion is Corollary 4.4 of \cite{mazurav}.  

Suppose now that $\ell \ne p$, $\L/\K$ is ramified, 
and $E$ has good reduction.  
Then $\K(E[p^\infty])/\K$ is unramified, so 
$\K(E(\L)[p^\infty]) = \K$, i.e., 
$
E(\K)[p^\infty] = E(\L)[p^\infty].
$
Now (ii) follows from Lemma \ref{parities}. 
\end{proof}

\begin{thm}
\label{par3}
Suppose that $v \nmid p$ and $E$ has good reduction at $v$.  
Let $w$ be a prime of $L$ above $v$.  
If $L_w/K_v$ is nontrivial and totally ramified, then 
$$\delta_v \equiv \dim_{\rf}E(K_v)[p] \pmod{2}.$$
\end{thm}

\begin{proof}
Let $L_w'$ be the intermediate field $K_v \subset L_w' \subset L_w$ 
with $[L_w:L_w'] = p$, as in Corollary \ref{intcor}.  
Applying Lemma \ref{newnormlem}(ii) to $L_w/L_w'$ and to $L_w'/K_v $ 
shows that
$$
N_{L_w/L_w'}E(L_w) = pE(L_w') \quad\text{and}\quad 
    E(K_v) \cap p E(L_w') = p E(K_v),
$$
so by Corollary \ref{intcor} and Lemma \ref{parities} we have 
$$
\delta_v \equiv \dim_{\rf}E(K_v)/p E(K_v) \equiv \dim_{\rf}E(K_v)[p] \pmod{2}.
$$
\end{proof}

\begin{thm}
\label{par5}
Suppose that $E$ is defined over $\Qp \subset K_v$
with good supersingular reduction at $p$.  
If $p = 3$ assume further that $|E(\F_3)| = 4$.

If $K_v$ contains 
the unramified quadratic extension of $\Qp$, then $\delta_v = 0$.
\end{thm}

\begin{proof}
Under these hypotheses $|E(\Fp)| = p+1$, so the characteristic polynomial 
of Frobenius on $E/\Fp$ is $X^2 + p$.  It follows that the 
characteristic polynomial of Frobenius over $E/\F_{p^2}$ is $(X+p)^2$.  In other 
words, multiplication by $-p$ reduces to the Frobenius endomorphism of 
$E/\F_{p^2}$

Let $\Qpp \subset K_v$ denote the unramified quadratic extension of $\Qp$, 
and $\Zpp$ its ring of integers.
Let $\Ehat$ denote the formal group over $\Zpp$ giving the kernel of 
reduction on $E$, and $[-p](X) \in \Zp[[X]]$ the power series giving 
multiplication by $-p$ on $\Ehat$.  Then $[-p](X) \equiv -pX \pmod{X^2}$, 
and since $-p$ reduces to Frobenius, we have $[-p](X) \equiv X^{p^2} \pmod{p}$.  
In other words, $\Ehat$ is a Lubin-Tate formal group of height $2$ over $\Zpp$, 
for the uniformizing parameter $-p$.

It follows that $\Zpp \subset \End(\Ehat)$.  Therefore $\Ehat(K_v)$ is a 
$\Zpp$-module, and since $E$ has supersingular reduction,  
$E(K_v)/pE(K_v) \cong \Ehat(K_v)/p\Ehat(K_v)$ is a vector space over 
$\Zpp/p\Zpp = \F_{p^2}$.  
Similarly, if $w$ is a prime of $L$ above $v$ then $\Ehat(L_w)$ is 
a $\Zpp[\Gal(L_w/K_v)]$-module and $E(L_w)/pE(L_w)$ is an $\F_{p^2}$-vector space.  
Hence  
$E(K_v)/(E(K_v) \cap N_{L_w/L_w'}E(L_w))$ is an $\F_{p^2}$-vector space, so 
its $\Fp$-dimension $\delta_v$ is even.
\end{proof}

\section{Dihedral extensions}
\label{diex}

Keep the notation of the previous sections.  
For cyclic extensions $L$ of $K$ in $F$, 
Proposition \ref{selfA} relates $\cork_{\Zp}\Sel_{p^\infty}(E/K)$ to 
$\dim_{\rf}\Sel_p(E/K)$, and 
Corollary \ref{rkdif1} relates $\dim_{\rf}\Sel_p(E/K)$ to $\dim_{\rf}\Sel_\P(A_L/K)$.   
Next we need to relate $\dim_{\rf}\Sel_\P(A_L/K)$ to 
$\cork_{\Zp}\Sel_{p^\infty}(A_L/K)$.
For this we need an additional hypothesis.

Suppose now that $c$ is an automorphism of order $2$ of $K$, let 
$\k \subset K$ be the fixed field of $c$, and suppose that $E$ is defined 
over $\k$.  Fix a cyclic extension $L/K$ of degree $p^n$, and let $A := A_L$, 
$R := R_L$, $\p \subset R$ the maximal ideal, 
etc., as in \S\ref{loc}.  We assume further that $L$ is Galois over $k$ with dihedral Galois 
group, i.e., $c$ acts by inversion on $G := \Gal(L/K)$.

\begin{thm}
\label{pol}
$
\dim_{\rf}(\Sh(A/K)/\Sh(A/K)_\div)[\P]
$
is even.
\end{thm}

Theorem \ref{pol} will be proved in Appendix \ref{polpf}.

\begin{rem}
\label{bad}
Theorem \ref{pol} is essential for our applications.  Without it, 
the formula in Proposition \ref{selfB2} below would not hold, and our 
approach would fail.  
The proof of Theorem \ref{pol} depends heavily on the 
fact that $L/\k$ is a dihedral extension.   
Stein \cite{stein} has given examples with $K = \Q$ where $L/\Q$ is abelian, 
$\Sh(A/\Q)$ is finite and $\dim_{\Fp}\Sh(A/\Q)[p]$ is odd.  

If $A$ had a polarization of degree prime to $p$, then 
Theorem \ref{pol} would follow directly from Tate's 
generalization of the Cassels pairing \cite{tate2}.  However, 
Howe \cite{howe} showed that (under mild hypotheses) 
every polarization of $A$ has degree divisible by $p^2$.
\end{rem}

Let $\Rp := R \otimes \Zp$.

\begin{prop}
\label{selfB2}
$$
\cork_{\Rp}\Sel_{p^\infty}(A/K) 
    \equiv \dim_{\rf}\HA(K,E[p]) - \dim_{\rf}E(K)[p] \pmod{2}.
$$
\end{prop}

\begin{proof}
The proof 
is identical to that of the formula for $\cork_{\Zp}\Sel_{p^\infty}(E/K)$ in 
Proposition \ref{selfA}, using Theorem \ref{athm}(ii) to view $R \subset \End_K(A)$, 
using Theorem \ref{pol} in place of the Cassels pairing, 
and using Proposition \ref{moreresprop} to identify $A(K)[\P]$ with $E(K)[p]$.
\end{proof}

\begin{thm}
\label{rkdif}
Suppose that $S$ is a set of primes of $K$ containing all primes above $p$, 
all primes ramified in $L/K$, and all primes where $E$ has bad 
reduction.  Then
$$
\cork_{\Zp}\Sel_{p^\infty}(E/K) - \cork_{\Rp}\Sel_{p^\infty}(A/K) 
   \equiv \sum_{v \in S} \delta_v \pmod{2}.
$$
\end{thm}

\begin{proof}
This follows directly from Corollary \ref{rkdif1} and Propositions 
\ref{selfA} and \ref{selfB2}.
\end{proof}

\begin{lem}
\label{unramlem}
Suppose $v$ is a prime of $K$ and $v = v^c$.  
Let $w$ be a prime of $L$ above $v$.  Then
\begin{enumerate}
\item
$L_w/K_v$ is totally ramified (we allow $L_w = K_v$),  
\item
if $v \nmid p$ and $L_w \ne K_v$ then $v$ is unramified in $K/\k$.
\end{enumerate}
\end{lem}

\begin{proof}
Let $w$ be a prime of $L$ above $v$, and $u$ the prime of $\k$ below $v$.  
Since $v = v^c$, the group $\Gal(L_w/\k_u)$ is dihedral.  
The inertia subgroup $I \subset \Gal(L_w/\k_u)$ 
is normal with cyclic quotient, and the only subgroups with this 
property are $\Gal(L_w/\k_u)$ and $\Gal(L_w/K_v)$.  This proves (i).  

Suppose now that $v$ is ramified in $K/\k$, and let $\ell$ be the 
residue characteristic of $K_v$.  By (i), 
the inertia group $I$ is a dihedral group of order $2[L_w:K_v]$.  
On the other hand, the Sylow $\ell$-subgroup of $I$ is normal with cyclic quotient (the 
tame inertia group).  The maximal abelian quotient of $I$ has order $2$, 
so $[L_w:K_v]$ must be a power of $\ell$, so $\ell = p$.
\end{proof}

\begin{lem}
\label{evenlem}
If $v$ is a prime of $K$ where $E$ has good reduction, $v \nmid p$, 
$v = v^c$, and $v$ is ramified in $L/K$, then 
$\dim_{\rf} E(K_v)[p]$ is even.
\end{lem}

\begin{proof}
Suppose $v \nmid p$, $v = v^c$, and $v$ ramifies in $L/K$.  
Fix a prime $w$ of $L$ above $v$, and let $u$ be the prime of $\k$ 
below $v$.  Let $\res_+$ and $\res$ denote the residue fields 
of $\k_u$ and $K_v$, respectively.  
Note that $K_v/\k_u$ is quadratic since $v = v^c$, 
and unramified by Lemma \ref{unramlem}(ii).  
Let $\phi$ be the Frobenius generator of $\Gal(K_v^\unr/\k_u)$, so 
$\phi^2$ is the Frobenius of $\Gal(K_v^\unr/K_v)$.

By Lemma \ref{unramlem}(i), $L_w/K_v$ is totally, tamely ramified.  
A standard result from algebraic number theory gives 
a $\Gal(\res/\res_+)$-equivariant injective homomorphism
$
\Gal(L_w/K_v) \hookto \res^\times.
$
Since $c$ acts by inversion on $\Gal(L_w/K_v)$, 
which is a nontrivial $p$-group by assumption, it follows that 
$\phi$ acts as inversion on $\bmu_p \subset \res^\times$.

Let $\alpha, \beta \in \bar\Fp^\times$ be the eigenvalues of  
$\phi$ acting on $E[p]$.  The Weil pairing and the action of $\phi$ on $\bmu_p$ 
show that $\alpha\beta = -1$.  
If $\alpha \ne \pm1$, then $1$ is not 
an eigenvalue of $\phi^2$ acting on $E[p]$, so 
$E(K_v)[p] = E[p]^{\phi^2 = 1} = 0$.
If $\alpha = \pm1$, then $\{\alpha,\beta\} = \{1,-1\}$, 
the action of $\phi$ on $E[p]$ is diagonalizable, $\phi^2$ 
is the identity on $E[p]$, and so $E(K_v)[p] = E[p]^{\phi^2 = 1} = E[p]$.
In either case, $\dim_{\rf}E(K_v)[p]$ is even.
\end{proof}

\begin{thm}
\label{par4}
If $v \mid p$ and $E$ has good ordinary reduction at $v$, 
then $\delta_v = 0$.
\end{thm}

\begin{proof}
Let $w$ be a prime of $L$ above $v$.  The theorem follows 
directly from Corollary \ref{intcor} and either 
Proposition \ref{oni} of Appendix \ref{appb} 
(if $L_w/K_v$ is totally ramified) or 
Lemma \ref{newnormlem}(i) (if not).
\end{proof}

\section{The main theorems}
\label{d1}

Fix a quadratic extension $K/\k$ with 
nontrivial automorphism $c$, an elliptic curve $E$ defined 
over $\k$, and an odd rational prime $p$.
Recall that if $F$ is an extension of $K$ then 
$\Scp(E/F) := \Hom(\Sel_{p^\infty}(E/F), \Qp/\Zp) \otimes \Qp$.
If $L$ is a cyclic extension of $K$ in $F$, let $R_L$ and $A_L$ 
be as defined in Definitions \ref{xirndef} and \ref{ildef}.

\begin{thm}
\label{second}
Suppose $F$ is an abelian $p$-extension of $K$, dihedral over $\k$ 
(i.e., $F$ is Galois over $\k$ and $c$ acts by inversion on $\Gal(F/K)$).  
Define
$$
\Sbad := \{\text{primes $v$ of $K$} : \text{$v$ ramifies in $F/K$ and $v = v^c$}\},
$$
and suppose that for every $v \in \Sbad$, one of the following three conditions holds:
\begin{enumerate}
\renewcommand{\theenumi}{(\alph{enumi})}
\item
$v \nmid p$ and $E$ has good reduction at $v$,
\item
$v \mid p$ and $E$ has good ordinary reduction at $v$,
\item
$v \mid p$, $E$ is defined over $\Qp \subset K_v$ 
with good supersingular reduction at $p$ 
(and if $p = 3$, then $|E(\F_3)| = 4$), and $K_v$ contains the unramified 
quadratic extension of $\Qp$.
\end{enumerate}
Then: 
\begin{enumerate}
\item
For every cyclic extension $L$ of $K$ in $F$, 
$$ 
\cork_{R_L \otimes \Zp}\Sel_{p^\infty}(A_L/K) 
    \equiv \cork_{\Zp}\Sel_{p^\infty}(E/K)  \pmod{2}.
$$
\item
If $\Xi$ is the set of cyclic extensions $L$ of $K$ contained in $F$, 
$G = \Gal(F/K)$, 
and $\Q[G] \cong \oplus_{L \in \Xi} \Q[G]_L$ is the decomposition \eqref{ssgd} 
of 
$\Q[G]$ into its isotypic components, then there an isomorphism of 
$\Qp[G]$-modules
$$
\Scp(E/F) \cong \dirsum{L \in \Xi} (\Q[G]_L \otimes \Qp)^{d_L}
$$
where for every $L$,
$$
d_L := \cork_{R_L \otimes \Zp}\Sel_{p^\infty}(A_L/K) 
    \equiv \cork_{\Zp}\Sel_{p^\infty}(E/K) \pmod{2}.
$$ 
\end{enumerate}
\end{thm}

\begin{proof}
Suppose that $L$ is a cyclic extension of $K$ in $F$, and 
let $R_p := R_L \otimes \Zp$ as in \S\ref{diex}.

Let $v$ be a prime of $K$.
If $v \ne v^c$ then $\delta_v + \delta_{v^c} \equiv 0 \pmod{2}$ by Lemma \ref{par2}.  
If $v = v^c$ and $v$ is unramified in $L/K$, then $v$ splits completely 
in $L/K$ by Lemma \ref{unramlem}(i), so $\delta_v = 0$ by 
Corollary \ref{intcor}.
Therefore by Theorem \ref{rkdif} we have 
$$
\cork_{\Zp}\Sel_{p^\infty}(E/K) - \cork_{R_p}\Sel_{p^\infty}(A_L/K) 
    \equiv \sum_{v \in \Sbad} \delta_v \pmod{2}.
$$
We will show that if $v \in \Sbad$ then $\delta_v = 0$, which will prove (i).

{\em Case 1: $v \nmid p$.}  
Then (a) holds, so $E$ has good reduction at $v$.  
If $w$ is a prime of $L$ above $v$, then $L_w/K_v$ is totally ramified 
by Lemma \ref{unramlem}(i).  Thus 
if $L_w = K_v$ then $\delta_v = 0$ by Corollary \ref{intcor}, and   
if $L_w \ne K_v$ then Theorem \ref{par3} and Lemma \ref{evenlem} show that 
$
\delta_v \equiv \dim_{\rf}E(K_v)[p] \equiv 0 \pmod{2}.
$

{\em Case 2: $v \mid p$.}  
Then either (b) or (c) must hold.
If (b) holds then $\delta_v = 0$ by Theorem \ref{par4}, and if (c) 
holds then $\delta_v = 0$ by Theorem \ref{par5}.  This proves (i).

By Corollary \ref{newressca}, 
$\Scp(E/F) \cong \oplus_{L \in \Xi} \;\Scp(A_L/K)$.  
By Theorem \ref{athm}(ii), 
$\Scp(A_L/K)$ is a vector space over the field 
$\Q[G]_L \otimes \Qp = R_L \otimes \Qp$,
and by (i) its dimension $d_L$ is congruent to 
$\cork_{\Zp}\Sel_{p^\infty}(E/K)$ modulo $2$.    
This proves (ii).
\end{proof}

\begin{thm}
\label{third}
Suppose $F/k$ and $E$ satisfy the hypotheses of Theorem \ref{second}.  

If $\cork_{\Zp}\Sel_{p^\infty}(E/K)$ is odd, then 
$\Scp(E/F)$ has a submodule isomorphic to 
$\Qp[\Gal(F/K)]$, and in particular 
$$
\cork_{\Zp}\Sel_{p^\infty}(E/F) \ge [F:K].
$$
\end{thm}

\begin{proof}
In Theorem \ref{second}(ii) we have $d_L \ge 1$ for every $L$, 
and the theorem follows.
\end{proof}

\begin{thm}
\label{fourth}
Suppose $F$ is an abelian $p$-extension of $K$, dihedral over $\k$, and 
all three of the following conditions are satisfied:
\begin{enumerate}
\renewcommand{\theenumi}{(\alph{enumi})}
\item
every prime $v \nmid p$ of $K$ that ramifies in $F/K$ satisfies $E(K_v)[p] = 0$,
\item
every prime $v$ of $K$ where $E$ has bad 
reduction splits completely in $F/K$,
\item
for every prime $v$ of $K$ dividing $p$, 
$E$ has good ordinary reduction at $v$ and if $\kappa$ is the 
residue field of $K_v$, then $E(\kappa)[p] = 0$.
\end{enumerate}
If $\Sel_{p^\infty}(E/K) \cong \Qp/\Zp$ 
(for example, if $\rk_\Z E(K) = 1$ and $\Sh(E/K)[p] = 0$), then 
$\Scp(E/F) \cong \Qp[\Gal(F/K)]$, 
and in particular $\cork_{\Zp}\Sel_{p^\infty}(E/F) = [F:K]$.
\end{thm}

\begin{proof}
Note that the hypotheses of this theorem are stronger than those 
of Theorem \ref{second}, so we can apply Theorem \ref{second}.

Suppose $L$ is a nontrivial cyclic extension of $K$ in $F$, 
and $K \subset M \subset L$ with $[L:M] = p$.  
We will show that for every prime $v$ of $K$ and $w$ of $L$ above $v$, 
\begin{equation}
\label{nc}
E(K_v) \subset N_{L_w/M_w} E(L_w).
\end{equation}
Assume this for the moment.  Then  
$\HA(K_v,E[p]) = \HE(K_v,E[p])$ for every $v$ by Proposition \ref{newint}, so 
if $\P_L$ is the prime above $p$ in $R_L \subset \End(A_L)$, 
we have
$$
\Sel_{\P_L}(A_L/K) = \HA(K,E[p]) = \HE(K,E[p]) = \Sel_p(E/K).
$$
Let $d_L := \cork_{R_L \otimes \Zp}\Sel_{p^\infty}(A_L/K)$.  
Using \eqref{selE} and \eqref{selEinf} (or the proof of 
Proposition \ref{selfA}) and Proposition \ref{moreresprop}, we have
$$
d_L \le \dim_{\Fp}\Sel_{\P_L}(A_L/K) - \dim_{\Fp}A_L[\P_L] = \dim_{\Fp}\Sel_p(E/K) 
    - \dim_{\Fp} E[p] = 1.
$$
But by Theorem \ref{second}(i), $d_L$ is odd, so $d_L = 1$.  
This holds for every $L$ (including $L = K$), so
the theorem follows directly from Theorem \ref{second}(ii).

It remains to prove \eqref{nc}.   

{\em Case 1:  $v \nmid p$, $E$ has good reduction at $v$, 
$v$ is unramified in $L/K$.} 
In this case \eqref{nc} holds by Lemma \ref{newnormlem}(i).

{\em Case 2: $v \nmid p$, $E$ has good reduction at $v$, 
$v$  is ramified in $L/K$.} 
In this case $E(K_v) = pE(K_v)$ 
by assumption (a) and Lemma \ref{parities}, so \eqref{nc} holds.

{\em Case 3:  $v \nmid p$, $E$ has bad reduction at $v$.} 
In this case $L_w = M_w$ by assumption (b), so 
\eqref{nc} holds.
 
{\em Case 4:  $v \mid p$.}  
If $L_w/K_v$ is not totally ramified, then $L_w/M_w$ is unramified 
and \eqref{nc} holds by Lemma \ref{newnormlem}(i).  If 
$L_w/K_v$ is totally ramified, then \eqref{nc} holds by 
Proposition \ref{oni} of Appendix \ref{appb} and assumption (c).
This completes the proof.
\end{proof}

\section{Special cases}
\label{examples}

\subsection{Odd Selmer corank}

In general it can be very difficult to determine the 
parity of $\cork_{\Zp}\Sel_{p^\infty}(E/K)$.  We now discuss some 
general situations in which the corank can be forced to be odd.

Fix an elliptic curve $E$ defined over $\Q$, and let $N_E$ be its 
conductor.  
Fix a Galois extension $K$ of $\Q$ such that $\Gal(K/\Q)$ 
is dihedral of order $2m$ with $m$ odd, $m \ge 1$.
Let $M$ be the quadratic extension of $\Q$ in $K$, $\Delta_M$ 
the discriminant of $M$, and $\chi_M$ 
the quadratic Dirichlet character attached to $M$.  
Let $c$ be one of the elements of order $2$ in $\Gal(K/\Q)$, 
and let $\k$ be the fixed field of $c$.  

\begin{lem}
\label{dclem}
$\cork_{\Zp}\Sel_{p^\infty}(E/K) \equiv \cork_{\Zp}\Sel_{p^\infty}(E/M) \pmod{2}$.
\end{lem}

\begin{proof}
The restriction map $\Scp(E/M) \to \Scp(E/K)^{\Gal(K/M)}$ is an isomorphism, so 
in the $\Qp$-representation $\Scp(E/K)/\Scp(E/M)$ of $\Gal(K/\Q)$, neither of 
the two one-dimensional representations occurs.  Since all other representations 
of $\Gal(K/\Q)$ have even dimension, we have that 
$$
\cork_{\Zp}\Sel_{p^\infty}(E/K) - \cork_{\Zp}\Sel_{p^\infty}(E/M) 
    = \dim_{\Qp}(\Scp(E/K)/\Scp(E/M))
$$ 
is even.
\end{proof}

The following proposition follows from the ``parity theorem'' for the 
$p$-power Selmer group proved by Nekov\'a\u{r} \cite{nekovarpc} and Kim \cite{kim}.

\begin{prop}
\label{selpar}
Suppose that $p > 3$ is a prime, and that $p$, $\Delta_M$, and $N_E$ 
are pairwise relatively prime.  Then 
$\cork_{\Zp}\Sel_{p^\infty}(E/K)$ is odd if and only if $\chi_M(-N_E) = -1$.
\end{prop}

\begin{proof}
Let $E'$ be the quadratic twist of $E$ by $\chi_M$, 
and let $w, w'$ be the signs in the functional equation of $L(E/\Q,s)$ 
and $L(E'/\Q,s)$, respectively.  Since $\Delta_M$ and $N_E$ 
are relatively prime, a well-known formula shows that $ww' = \chi_M(-N_E)$.  

Using Lemma \ref{dclem} we have 
\begin{align*}
\cork_{\Zp}\Sel_{p^\infty}(E/K) &\equiv \cork_{\Zp}\Sel_{p^\infty}(E/M) \pmod{2}\\
    &= \cork_{\Zp}\Sel_{p^\infty}(E/\Q) + \cork_{\Zp}\Sel_{p^\infty}(E'/\Q).
\end{align*}
By a theorem of Nekov\'a\u{r} \cite{nekovarpc} (if $E$ has ordinary reduction at $p$) 
or Kim \cite{kim} (if $E$ has supersingular reduction at $p$), 
we have that $\cork_{\Zp}\Sel_{p^\infty}(E/\Q)$ is even if and only if 
$w = 1$, and similarly for $E'$ and $w'$.  Thus $\cork_{\Zp}\Sel_{p^\infty}(E/K)$ 
is odd if and only if $w = -w'$, and the proposition follows.
\end{proof}

For every prime $p$, let $\Kantip$ be the maximal 
abelian $p$-extension of $K$ that is Galois and dihedral 
over $\k$, and unramified (over $K$) at all primes dividing $N_E$ that 
do not split in $M/\Q$.
(Note that if a rational prime $\ell$ splits in $M$, then every 
prime of $\k$ above $\ell$ splits in $K/\k$ since $[K:M]$ is odd.)

\begin{thm}
\label{dihex}
Suppose $p > 3$ is prime, and $p$, $\Delta_M$, and $N_E$ 
are pairwise relatively prime.  
If $\chi_M(-N_E) = -1$, 
then for every finite extension $F$ of $K$ in $\Kantip$, 
$$
\cork_{\Zp}\Sel_{p^\infty}(E/F) \ge [F:K].
$$
\end{thm}

\begin{proof}
This will follow directly from Theorem \ref{third} and Proposition \ref{selpar}, 
once we show that the hypotheses of Theorem \ref{second} are satisfied.   
By definition of $\Kantip$, the set $\Sbad$ of Theorem \ref{second} 
contains only primes above $p$, and since $p \nmid N_E \Delta_M$ either 
(b) or (c) holds for every $v \in \Sbad$.
\end{proof}

If $m = 1$, so $K = M$, and if $M$ is imaginary, then $\Kantip$ contains 
the anticyclotomic $\Zp$-extension of $K$, and thanks to \cite{cornut, vatsal} 
we know that the bulk of the contribution to the Selmer groups in 
Theorem \ref{dihex} comes from Heegner points.

If $m = 1$ and $M$ is real, then there is no $\Zp$-extension of $K$ 
in $\Kantip$.  However, $\Kantip$ is still an infinite extension of $K$, 
and (for example) 
every finite abelian $p$-group occurs as a quotient of $\Gal(\Kantip/K)$.  

More generally, for arbitrary $m$, if $M$ is imaginary then $\Kantip$ contains 
a $\Zp^d$-exten\-sion of $K$ with $d = (m+1)/2$, and if $M$ is real then 
$K$ is totally real so $\Kantip$ is 
infinite but contains no $\Zp$-extension of $K$.  Except for Heegner points 
in special cases (such as when $m = 1$ and $M$ is imaginary), it is not known where 
the Selmer classes in Theorem \ref{dihex} come from.

\subsection{Split multiplicative reduction at $p$}
Suppose now that $K/\k$ is a quadratic extension, and $F$ is a finite abelian 
$p$-extension of $K$, dihedral over $\k$.  Suppose that $E$ is an elliptic 
curve over $\k$, and $v$ is a prime of $K$ above $p$, inert in $K/\k$, where $E$ has 
split multiplicative reduction.  If $F/K$ is ramified at $v$ then 
Theorems \ref{second} and \ref{third} do not apply.  We now study this case more 
carefully.

\begin{lem}
\label{badlem}
Suppose $v$ is a prime of $K$ above $p$ such that $v = v^c$, $u$ is the 
prime of $\k$ below $v$, and  
$E$ has split multiplicative reduction at $u$.  If $L$ is a nontrivial cyclic 
extension of $K$ in  $F$, $v$ is totally ramified in $L/K$, 
$K \subset L' \subset L$ with $[L:L'] = p$, 
and $w$ is a prime of $L$ above $v$, then 
$
[E(K_v) : E(K_v) \cap N_{L/L'}E(L_w)] = p.
$
\end{lem} 

\begin{proof}
Let $\m_{u}$ denote the maximal ideal of $\k_u$.
Since $E$ has split multiplicative reduction, there is a nonzero $q \in \m_{u}$ 
such that $E(L_w) \cong L_w^\times / q^\Z$ as $\Gal(L_w/\k_u)$-modules.  

Since $v = v^c$, $L_w/\k_v$ is dihedral so the maximal 
abelian extension of $\k_v$ in $L_w$ is $K_v$.  
Thus local class field theory gives an identity of norm groups
$$
N_{K_v/\k_v}K_v^\times = N_{L_w/\k_v}L_w^\times 
    \subset N_{L_w/L_w'}L_w^\times.
$$
Since $q^2 \in N_{K_v/\k_v}K_v^\times$ and 
$[(L_w')^\times:N_{L_w/L_w'}L_w^\times] = [L_w:L_w'] =p$ is odd, 
we see that $q \in N_{L_w/L_w'}L_w^\times$, and so
\begin{equation}
\label{normin}
[E(K_v) : E(K_v) \cap N_{L/L'}E(L_w)] 
    = [K_v^\times : K_v^\times \cap N_{L/L'} L_w^\times].
\end{equation}
Let $[L:K] = p^n$.  
If $\spair{\;}{\;}$ denotes the Artin map of local class field theory, 
then $K_v^\times \cap N_{L/L'} L_w^\times$ is the kernel of the 
map $K_v^\times \to \Gal(L_w/K_v)$ given by 
$$
x \mapsto [x,L_w/L'_w] = [N_{L'/K}x,L_w/K_v] = [x^{p^{n-1}},L_w/K_v] 
    = [x,L_w/K_v]^{p^{n-1}}.
$$
Since $x \mapsto [x,L_w/K_v]$ maps $K_v^\times$ onto a cyclic group of 
order $p^n$, we conclude that the index \eqref{normin} is $p$, as desired.
\end{proof}

Let $\Sbad_p$ be the set of primes $v$ of $K$ above $p$ such that $v = v^c$ and  
neither of the hypotheses (b) or (c) of Theorem \ref{second} hold for $v$.

\begin{thm}
\label{smult}
Suppose that $F$ is a finite abelian $p$-extension of $K$ that is dihedral over $\k$ 
and unramified at all primes $v \nmid p$ of bad reduction that do not split in 
$K/\k$.  Suppose further that for every prime $v \in \Sbad_p$, 
$E$ has split multiplicative reduction at $v$ and $v$ is totally ramified in $F/K$.  
Then:
\begin{enumerate}
\item
If $\cork_{\Zp}\Sel_{p^\infty}(E/K) + |\Sbad_p|$ is odd, then 
$$
\cork_{\Zp}\Sel_{p^\infty}(E/F) \ge \cork_{\Zp}\Sel_{p^\infty}(E/K) + [F:K] - 1.
$$
\item
If $\Sel_{p^\infty}(E/K)$ is finite and $|\Sbad_p|$ is odd, then 
$$
\cork_{\Zp}\Sel_{p^\infty}(E/F) \ge [F:K] - 1.
$$
\item
Suppose that $|\Sbad_p| = 1$, and the hypotheses (a), (b), or (c) 
of Theorem \ref{fourth} hold for every prime $v$ of $K$ not in $\Sbad_p$. 
If $\Sel_{p^\infty}(E/K) = 0$, then 
$$
\cork_{\Zp}\Sel_{p^\infty}(E/F) = [F:K] - 1.
$$
\end{enumerate}
\end{thm}

\begin{proof}
The proof is identical to that of Theorems \ref{third} and \ref{fourth}, except that we 
use Lemma \ref{badlem} to compute the $\delta_v$ for $v \in \Sbad_p$.

Suppose $L$ is a nontrivial cyclic extension of $K$ in $F$.  
Exactly as in Theorem \ref{second}, we have 
$\sum_{v \notin \Sbad_p} \delta_v \equiv 0 \pmod{2}$.
If $v \in \Sbad_p$, then $\delta_v = 1$ by Lemma \ref{badlem} and 
Corollary \ref{intcor}.  Thus we conclude that 
$\sum_v \delta_v \equiv |\Sbad_p| \pmod{2}$.
Exactly as in Theorem \ref{second} we conclude using Theorem \ref{rkdif} that 
\begin{equation}
\label{rl}
\Scp(E/F) \cong \dirsum{L \in \Xi} (\Q[G]_L \otimes \Qp)^{d_L}
\end{equation}
where $d_L \equiv \cork_{\Zp}\Sel_{p^\infty}(E/K) + |\Sbad_p| \pmod{2}$ 
for every $L \ne K$.  
Assertion (i) now follows 
exactly as in the proof of Theorem \ref{third}, and (ii) is a special case of (i).

For (iii), it follows exactly as in the proof of Theorem \ref{fourth} that 
$\HA(K_v,E[p]) = \HE(K_v,E[p])$ for every $v \notin \Sbad_p$.  
Thus if $\Sbad_p = \{v_0\}$, there is an exact sequence
\begin{equation}
\label{exse}
0 \to \Hs{\E\cap\A}(K,E[p]) \to \HA(K,E[p]) 
    \to \HA(K_{v_0},E[p])/\Hs{\E\cap\A}(K_{v_0},E[p]).
\end{equation}
By Lemma \ref{badlem} and Proposition \ref{newint},
$$
\dim_{\Fp}\HA(K_{v_0},E[p]) = \dim_{\Fp}\HE(K_{v_0},E[p]) 
    = \dim_{\Fp}\Hs{\E\cap\A}(K_{v_0},E[p]) + 1
$$
(the first equality holds because $\A$ and $\E$ are self-dual), 
so it follows from \eqref{exse} that 
\begin{multline*}
\dim_{\Fp} \Sel_{\P_L}(A_L/K) = \dim_{\Fp}\HA(K,E[p]) 
    \le \dim_{\Fp}\HE(K,E[p]) + 1 \\
    = \dim_{\Fp}E[p] + 1 = \dim_{\Fp}A[\P] + 1.
\end{multline*}
Therefore $d_L := \cork_{R_L \otimes \Zp}\Sel_{p^\infty}(A_L/K) \le 1$.  
The proof of (i) showed that $d_L$ is odd, so $d_L = 1$.  
Hence in \eqref{rl} we have $d_L = 1$ if $L \ne K$, and $d_K = 0$.  This proves (iii).
\end{proof}

\begin{rem}
In the case where $K = M$ is imaginary quadratic and $F$ is 
a subfield of the anticyclotomic $\Zp$-extension, Bertolini and Darmon 
\cite{bd} give a construction of Heegner-type points that account for 
most of the Selmer classes in Theorem \ref{smult}.
\end{rem}

\appendix

\section{Skew-Hermitian pairings}
\label{la}
\label{polpf}

In this appendix we prove Proposition \ref{selfB} and Theorem \ref{pol}.

Let $p$ be an odd prime, 
$L/K$ be a cyclic extension of number fields of degree $p^n$, 
$G := \Gal(L/K)$, and $\cR := R_L \otimes \Zp$, where $R_L$ is given 
by Definition \ref{xirndef}.  
We view $\cR$ as a $G_K$-module by letting $G_K$ act 
{\em trivially} (not the action induced from the action on $R_L$).  
Then $\cR$ is the cyclotomic ring over $\Zp$ generated by $p^n$-th 
roots of unity (see for example \cite{prim} Lemma 5.4(ii)).  

Let $\iota$ be the involution of $R_L$ (resp., $\cR$) 
induced by $\zeta \mapsto \zeta^{-1}$ for $p^n$-th roots of unity 
$\zeta \in R_L$ (resp., $\zeta \in \cR$).
If $\T$ is an $\cR$-module, we let $\T^\iota$ be the $\cR$-module 
whose underlying abelian group is $\T$, but with $\cR$-action 
twisted by $\iota$.  

\begin{defn}
\label{shdef}
Suppose $\T$ is an $\cR$-module and $B$ is a $\Zp$-module.  We say that 
a $\Zp$-bilinear pairing  
$$\pair{\;}{\;} : \T \times \T \to B$$ 
is {\em $\iota$-adjoint} 
if $\pair{rx}{y} = \pair{x}{r^\iota y}$ for every $r \in \cR$ and $x,y \in \T$.
We say that a pairing  
$$\pair{\;}{\;} : \T \times \T \to \cR \otimes_{\Zp} B$$ 
is {\em $\cR$-semilinear} if 
$\pair{rx}{y} = r\pair{x}{y} = \pair{x}{r^\iota y}$ for every 
$r \in \cR$ and $x, y \in \T$,
and we say $\pair{\;}{\;}$ is {\em skew-Hermitian} if it is $\cR$-semilinear and 
$\pair{y}{x} = -\pair{x}{y}^{\iota \otimes 1}$ for every $x, y \in \T$.

We say that $\pair{\;}{\;}$ is {\em nondegenerate} (resp., {\em perfect}) 
if the induced map $\T \to \Hom_{\Zp}(\T^\iota,B)$ 
(or $\Hom_\cR(\T^\iota,\cR \otimes_{\Zp} B)$, 
depending on the context) is injective (resp., an isomorphism).
\end{defn}

\begin{defn}
\label{pdefs}
Let $\zeta$ be a primitive $p^n$-th root of unity in 
$R_L$, and let $\pi := \zeta-\zeta^{-1}$.  Then $\pi$ is a generator 
of the prime $\P_L$ of $R_L$ above $p$, and $\pi$ is also a generator of 
the maximal ideal $\P$ of $\cR$, and $\pi^\iota = -\pi$.
Let $d := \pi^{p^{n-1}(pn-n-1)}$, so $d$ is a generator of 
the inverse different of $R_L/\Z$ and of $\cR/\Zp$, and $d^\iota = -d$.  
Define a trace pairing
$$
t_{\cR/\Zp} : \cR \times \cR \to \Zp, 
    \quad t_{\cR/\Zp}(r,s) := \Tr_{\cR/\Zp}(d^{-1} r s^\iota)
$$
This pairing is $\iota$-adjoint, perfect, and (since $d^\iota = -d$) 
skew-symmetric.  
Define $\tau : \cR \to \Zp$ by 
$
\tau(r) := t_{\cR/\Zp}(1,r) = -\Tr_{\cR/\Zp}(d^{-1}r).
$
\end{defn}

\begin{lem}
\label{g1}
Suppose that $\T$ is an $\cR[G_K]$-module and $B$ is a $\Zp[G_K]$-module.  
Composition with $\tau \otimes 1 : \cR \otimes_{\Zp} B \to B$ 
gives an isomorphism of $G_K$-modules
$$\Hom_{\cR}(\T, \cR \otimes_{\Zp} B) \isom \Hom_{\Zp}(\T,B).$$
\end{lem}

\begin{proof}
We will construct an inverse to the map in the statement of the lemma.  
Suppose $f \in \Hom_{\Zp}(\T,B)$.  Fix a $\Zp$-basis $\{\nu_1,\ldots,\nu_b\}$ 
of $\cR$, and let $\{\nu_1^*,\ldots,\nu_b^*\}$ be the dual basis with 
respect to $t_{\cR/\Zp}$, i.e., $t_{\cR/\Zp}({\nu_i},{\nu_j^*}) = \delta_{ij}$.
For $x \in \T$ define
$$
\hat{f}(x) := \sum_{i=1}^b \nu_i^* \otimes f(\nu_i^\iota x) \in \cR \otimes_{\Zp} B.
$$
Then for every $j$ and $x$,
$$
(\tau \otimes 1)(\nu_j^\iota \hat{f}(x)) 
    = \sum_{i=1}^b t_{\cR/\Zp}({1},{\nu_j^\iota\nu_i^*}) f(\nu_i^\iota x) 
    = \sum_{i=1}^b t_{\cR/\Zp}({\nu_j},{\nu_i^*}) f(\nu_i^\iota x) = f(\nu_j^\iota x)
$$
Since the $\nu_j$ are a basis of $\cR$, we conclude that 
\begin{equation}
\label{iofb}
\text{$(\tau \otimes 1)(r \hat{f}(x)) = f(r x)$ for every $r \in \cR$}.
\end{equation}
Thus if $s \in \cR$ then for every $r$
$$
(\tau \otimes 1)(r \hat{f}(sx)) = f(rsx) = (\tau \otimes 1)(rs \hat{f}(x)).
$$
Since $t_{\cR/\Zp}$ is perfect and $\cR$ is free over $\Zp$, 
it follows that $\hat{f}(sx) = s\hat{f}(x)$, so 
$\hat{f} \in \Hom_{\cR}(\T,\cR \otimes_{\Zp} B)$.

By \eqref{iofb} with $r=1$, $(\tau \otimes 1) \circ \hat{f} = f$ so 
$\Hom_{\cR}(\T,\cR \otimes_{\Zp} B) \map{\circ(\tau\otimes 1)} \Hom_{\Zp}(\T,B)$ 
is surjective.  The injectivity follows from the fact that 
$t_{\cR/\Zp}$ is perfect and $\cR$ is free over $\Zp$.  
The $G_K$-equivariance is clear (recall that $G_K$ acts trivially on $\cR$).
\end{proof}

\begin{prop}
\label{gorprop}
Suppose that $\T$ is an $\cR$-module and $B$ is a $\Zp$-module.  
Composition with $\tau \otimes 1 : \cR \otimes_{\Zp} B \to B$ gives a 
bijection between the set of $\cR$-semilinear pairings 
$\T \times \T \to \cR \otimes_{\Zp} B$, and the set of 
$\iota$-adjoint pairings 
$\T \times \T \to B$.

If $\pair{\;}{\;}_\cR$ maps to $\pair{\;}{\;}_{\Zp}$ under this bijection, 
then $\pair{\;}{\;}_{\Zp}$ is perfect (resp., $G_K$-equivariant) 
if and only if $\pair{\;}{\;}_\cR$ is perfect (resp., $G_K$-equivariant).
\end{prop}

\begin{proof}
By Lemma \ref{g1}, composition with $\tau \otimes 1$ induces 
a $G_K$-isomorphism
\begin{equation}
\label{dps}
\Hom_{\cR}(\T,\Hom_{\cR}(\T^\iota,\cR \otimes_{\Zp} B)) \isom 
    \Hom_{\cR}(\T,\Hom_{\Zp}(\T^\iota,B)).
\end{equation}
The left-hand side is the set of 
$\cR$-semilinear pairings $\T \times \T \to \cR \otimes_{\Zp} B$,
and the right-hand side is the set of $\iota$-adjoint pairings 
$\T \times \T \to B$.

Since composition with $\tau \otimes 1$ identifies 
the isomorphisms in \eqref{dps}, we see that $\pair{\;}{\;}_\cR$ 
is perfect if and only if $\pair{\;}{\;}_{\Zp}$ is perfect.  Since 
\eqref{dps} is $G_K$-equivariant, $\pair{\;}{\;}_\cR$ 
is $G_K$-equivariant if and only if $\pair{\;}{\;}_{\Zp}$ is $G_K$-equivariant.
This completes the proof of the proposition.
\end{proof}

Let $A$ be the abelian variety $A_L$ of Definition \ref{ildef}.
Recall (Definitions \ref{pdefs} and \ref{xirndef} and Theorem \ref{athm}(i)) 
that $\pi$ is a generator of the prime $\P_L$ of $R_L$, $\pi\cR = \P$, 
and $\I_L = \P_L^{p^{n-1}}$.  

\begin{defn}
\label{lastpair}
Define a pairing $f : \I_L \times \I_L \to R_L$ by 
$$
f(\alpha,\beta) := \pi^{-2p^{n-1}} \alpha \beta^\iota.
$$
Theorem \ref{athm}(iv) gives a $G_{K}$-isomorphism 
$T_p(A) \cong \I_L \otimes T_p(E)$, and using this identification we define
$$
\pair{\;}{\;}_\cR := f \otimes e : T_p(A) \times T_p(A) \too \cR \otimes_{\Zp} \Zp(1)
$$ 
where $e$ is the Weil pairing on $E$.  
In other words, if $\alpha, \beta \in \I_L$ and $x,y \in T_p(E)$, we set
$$
\pair{\alpha \otimes x}{\beta \otimes y}_\cR 
    := (\pi^{-2p^{n-1}} \alpha \beta^\iota) \otimes e(x,y).
$$
\end{defn}

\begin{lem}
\label{lastlem}
The pairing $\pair{\;}{\;}_\cR$ of Definition \ref{lastpair} is perfect, 
$G_K$-equivariant, and skew-Hermitian.
\end{lem}

\begin{proof}
The Weil pairing is perfect and skew-symmetric, and   
the pairing $f$ is perfect and Hermitian (since $\pi^\iota = -\pi$).  
Thus $\pair{\;}{\;}_\cR$ is perfect and skew-Hermitian.
If $\alpha,\beta \in \I_L$, $x,y \in T_p(E)$, and $\gamma \in G_K$ then
\begin{align*}
\pair{(\alpha \otimes x)^\gamma}{(\beta \otimes y)^\gamma}_\cR 
    &= \pair{\alpha \gamma^{-1} \otimes \gamma x}{\beta \gamma^{-1} \otimes \gamma y}_\cR \\
    &= \pi^{-2p^{n-1}} (\alpha \gamma^{-1})(\beta \gamma^{-1})^\iota 
        \otimes e(\gamma x, \gamma y) \\
    &= \pi^{-2p^{n-1}} (\alpha \gamma^{-1})(\beta^\iota \gamma) \otimes e(x,y)^\gamma \\
    &= f(\alpha,\beta) \otimes e(x,y)^\gamma 
    = \pair{\alpha \otimes x}{\beta \otimes y}_\cR^\gamma
\end{align*}
since the Weil pairing is $G_K$-equivariant and $G_K$ acts trivially on $\cR$.  
\end{proof}

The following is Proposition \ref{selfB}.

\begin{prop}
\label{sd}
The Selmer structure $\A$ on $E[p]$ of Definition \ref{bdef} is 
self-dual.
\end{prop}

\begin{proof}
Using Proposition \ref{gorprop}, we let 
$$
\pair{\;}{\;}_{\Zp} :  T_p(A) \times T_p(A) \too \Zp(1)
$$
be the pairing corresponding under Proposition \ref{gorprop} to 
the pairing $\pair{\;}{\;}_\cR$ of Definition \ref{lastpair}, with $B = \Zp(1)$.  
It follows from Proposition \ref{gorprop} and Lemma \ref{lastlem} 
that $\pair{\;}{\;}_{\Zp}$ 
is perfect, $G_K$-equivariant, and $\iota$-adjoint.

By a generalization of Tate duality due to Bloch and Kato 
(see Proposition 3.8 and Example 3.11 of \cite{bk}), 
for every prime $v$ of $K$, the  
pairing $\pair{\;}{\;}_{\Zp}$ induces a perfect, $\iota$-adjoint cup-product pairing 
$$
\lambda : H^1(K_v,T_p(A)) \times H^1(K_v,T_p(A) \otimes \Qp/\Zp) \too \Qp/\Zp.
$$
and under this pairing the image of 
$A(K_v) \to H^1(K_v,T_p(A))$ and the image of 
$A(K_v) \otimes \Qp/\Zp \to H^1(K_v,T_p(A) \otimes \Qp/\Zp)$ 
are orthogonal complements of each other.  

The pairing $\lambda$ induces a pairing 
$$
\lambda_{\P_L} :  H^1(K_v,T_p(A)/\P_L T_p(A)) 
    \times H^1(K_v,(T_p(A) \otimes \Qp/\Zp)[\P_L]) \too \Fp.
$$
We have isomorphisms (the first one uses the chosen generator $\pi$ of $\P_L$)
$$
T_p(A)/\P_L T_p(A) \cong A[\P_L] \cong (T_p(A) \otimes \Qp/\Zp)[\P_L].
$$
Along with the identification $A[\P_L] \cong E[p]$ of Proposition \ref{moreresprop},  
this transforms $\lambda_{\P_L}$ into a pairing 
$H^1(K_v,E[p]) \times H^1(K_v,E[p]) \to \Fp$, 
and one can check directly from the definition of $\lambda$ that 
this pairing is the same as the local cup product pairing on $H^1(K_v,E[p])$ 
coming from the Weil pairing as in \S\ref{var} and \S\ref{ecx}.

A couple of straightforward diagram chases (see for example 
Lemma 1.3.8 and Proposition 1.4.3 of \cite{EulerSystems}) 
show that the image of 
\begin{equation}
\label{dsnt}
A(K_v) \to H^1(K_v,T_p(A)) \to H^1(K_v,T_p(A)/\P_L T_p(A)) \isom H^1(K_v,E[p])
\end{equation}
and the inverse image of $A(K_v) \otimes \Qp/\Zp$ under 
$$
H^1(K_v,E[p]) \isom H^1(K_v,(T_p(A) \otimes \Qp/\Zp)[\P_L]) 
    \to H^1(K_v,T_p(A) \otimes \Qp/\Zp)
$$
are equal and are orthogonal complements under $\lambda_{\P_L}$.  
By definition the image of \eqref{dsnt} is $\HA(K_v,E[p])$, 
so this proves that $\A$ is self-dual.
\end{proof}

It remains to prove Theorem \ref{pol}, and for that we need to 
be in the dihedral setting of \S\ref{diex}.  We assume now that $K$ 
has an automorphism $c$ of order $2$, that $E$ is defined over the 
fixed field $\k$ of $K$, that $L$ is Galois over $\k$, and that $c$ acts 
by inversion on $G := \Gal(L/K)$.

We begin by fixing a model of $A$ defined over $\k$.  

\begin{defn}
\label{newadef}
Fix a lift of $c$ to $G_{\k}$, and denote this lift by $\c$.  
Then $\Gal(L/\k)$ is the semidirect product $G \semidirect H$, 
where $H$ is the group of order $2$ generated by the restriction 
of $\c$.  Let $\J_L := (1+\c)\I_L$, where 
$\I_L \subset \Z[G] \subset \Z[\Gal(L/\k)]$ is 
the ideal of $\Z[G]$ given in Definition \ref{xirndef}.  
Then $\J_L$ is a right ideal of $\Z[\Gal(L/\k)]$, 
and we define an abelian variety $A'$ over $\k$ by
$$
A' := \J_L \otimes E
$$
as in Definition 1.1 (and \S6) of \cite{prim}. 
\end{defn}

\begin{prop}
\label{secondtolast}
\begin{enumerate}
\item
Left multiplication by $(1+\c)$ is an isomorphism of right $G_K$-modules  
from $\I_L$ to $\J_L$.
\item
The isomorphism of (i) induces an isomorphism 
$A \cong A'$ defined over $K$.  
\end{enumerate}
\end{prop}

\begin{proof}
The first assertion is easily checked, and the second 
follows by Corollary 1.9 of \cite{prim}.  See also Theorem 6.3 of \cite{prim}.
\end{proof}

From now on we view $A$ as defined over $\k$, by using the model 
$A'$ of $A$ and Proposition \ref{secondtolast}(ii).  
We extend the $G_K$-action on $T_p(A)$, $\I_L$, and $\cR$ to a $G_\k$-action 
by identifying $T_p(A)$ with $T_p(A')$, $\I_L$ with $\J_L$ as in 
Proposition \ref{secondtolast}(i), and 
letting $\c$ act on (the trivial $G_K$-module) $\cR$ by $\iota$.  
The actions on $T_p(A)$ and $\I_L$ depend on the choice of $\c$.

\begin{prop}
\label{lastprop}
With the conventions above, the pairing 
$$\pair{\;}{\;}_\cR : T_p(A) \times T_p(A) \to \cR \otimes_{\Zp} \Zp(1)$$ 
of Definition \ref{lastpair} is $G_\k$-equivariant.
\end{prop}

\begin{proof}
By Theorem 2.2(iii) of \cite{prim}, 
there is a $G_\k$-isomorphism $T_p(A') \cong \J_L \otimes T_p(E)$. 
With the conventions above, this says that the isomorphism 
$T_p(A) \cong \I_L \otimes T_p(E)$, which was used to construct the pairing 
$\pair{\;}{\;}_\cR$ in Definition \ref{lastpair}, is $G_\k$-equivariant.
The proposition follows from this exactly as in Lemma \ref{lastlem}, 
using the fact that for $\alpha \in \I_L$, 
$(1+\c)\alpha\c = (1+\c)\c\alpha^\iota = (1+\c)\alpha^\iota$.
\end{proof}

Let $\Dp := \cR \otimes_{\Zp} \Qp/\Zp$.  

\begin{prop}
\label{linalg}
Suppose that $\T$ is an $\cR$-module of finite cardinality 
and a $\Gal(K/\k)$-module, 
and suppose that there is a nondegenerate, skew-Hermitian, 
$\Gal(K/\k)$-equivariant pairing 
$$
\spair{\;}{\;} : \T \times \T \too \Dp.
$$  
Then $\T$ has isotropic $\cR$-submodules $M$, $M'$ such that $M \cong M'$ and  
$
\T = M \oplus M'.
$
In particular $\dim_{\Fp}\T[\p]$ is even.
\end{prop}

\begin{proof}
Define a pairing $\spair{\;}{\;}' : \T \times \T \to \Dp$ by 
$\spair{v}{w}' := \spair{v}{c w}$.  It is straightforward to check 
that the pairing $\spair{\;}{\;}'$ is non-degenerate, $\cR$-bilinear, 
and skew-symmetric.  The proposition now follows by a well-known argument.
\end{proof}

We will now deduce Theorem \ref{pol} from a  
(slight generalization of a) result of Flach.
Let $\Shmd := \Sh(A/K)[p^\infty]/\Sh(A/K)[p^\infty]_\div$.  

\begin{thm}[Flach \cite{flach}]
\label{flachthm}
Suppose that 
$$
\bpair{\;}{\;}_\cR : T_p(A) \times T_p(A) \to \cR \otimes_{\Zp} \Zp(1)
$$ 
is a perfect, $G_{\k}$-equivariant, skew-Hermitian pairing.  
Then there is a perfect, $\Gal(K/\k)$-equivariant, skew-Hermitian pairing, 
$$
\spair{\;}{\;}_\cR : \Shmd \times \Shmd \to \Dp.
$$
\end{thm}

\begin{proof}
This is essentially Theorems 1 and 2 of \cite{flach}.  
We sketch here the minor modifications to the arguments of \cite{flach} needed to 
prove Theorem \ref{flachthm}.

Given a $G_K$-equivariant pairing 
$T_p(A) \times T_p(A) \to \Zp(1)$, Flach constructs a pairing 
$\Shmd \times \Shmd \to \Qp/\Zp$.  The definition (\cite{flach} p.\ 116) 
is given explicitly in terms of cocycles.  Since 
$G_K$ acts trivially on $\cR$, we have canonical isomorphisms 
\begin{equation}
\label{func}
H^i(K, \cR \otimes_{\Zp} \Zp(1)) \cong \cR \otimes_{\Zp} H^i(K,\Zp(1))
\end{equation}
for every $i$, and similarly with $K$ replaced by any of its completions $K_v$ 
and/or with $\Zp(1)$ replaced by $\Qp/\Zp(1)$.  
The isomorphisms \eqref{func} come from analogous isomorphisms on modules 
of cocycles.
Using this, starting with our pairing $\bpair{\;}{\;}_\cR$ and 
following Flach's construction verbatim produces a pairing 
$$
\spair{\;}{\;}_\cR : \Shmd \times \Shmd \to \Dp.
$$  
We need to show that 
$\spair{\;}{\;}_\cR$ is perfect, $\Gal(K/\k)$-equivariant, and skew-Hermitian.

The fact that $\spair{\;}{\;}_\cR$ is $\Gal(K/\k)$-equivariant 
follows directly from the definition in \cite{flach}, 
as each step is canonical and Galois-equivariant.

Similarly, following the definition in \cite{flach} and using that 
$\bpair{\;}{\;}_\cR$ is skew-Hermitian, one sees directly 
that $\spair{rx}{y}_\cR = r\spair{x}{y}_\cR = \spair{x}{r^\iota y}_\cR$ 
for every $r \in R$, $x,y \in \Shmd$. 
The fact that $\spair{y}{x}_\cR = -\spair{x}{y}_\cR^\iota$ is proved exactly as 
Theorem 2 of \cite{flach}, which proves the skew-symmetry of the pairing 
in Flach's setting.

It remains only to show that $\spair{\;}{\;}_\cR$ is perfect, or equivalently 
(since $\Shmd$ is finite) $\spair{\;}{\;}_\cR$ is nondegenerate.  
Let $\bpair{\;}{\;}_{\Zp} : T_p(A) \times T_p(A) \to \Zp(1)$ 
(resp., $\spair{\;}{\;}_{\Zp} : \Shmd \times \Shmd \to \Qp/\Zp$) be the pairing 
corresponding to $\bpair{\;}{\;}_\cR$ (resp., $\spair{\;}{\;}_\cR$) 
under the correspondence of Proposition \ref{gorprop}.  

By Proposition \ref{gorprop}, since $\bpair{\;}{\;}_\cR$ is perfect, 
$\bpair{\;}{\;}_{\Zp}$ is perfect.  
One can check from the definition that $\spair{\;}{\;}_{\Zp}$ is the pairing 
Flach constructs from $\bpair{\;}{\;}_{\Zp}$, and thus Flach's Theorem 1 
shows that $\spair{\;}{\;}_{\Zp}$ is perfect.  Now Proposition \ref{gorprop}
shows that $\spair{\;}{\;}_\cR$ is perfect.
This completes the proof of the theorem.
\end{proof}

\begin{proof}[Proof of Theorem \ref{pol}]
We apply Theorem \ref{flachthm}, 
using the pairing $\pair{\;}{\;}_\cR$ of Definition \ref{lastpair} 
(along with Lemma \ref{lastlem} and Proposition \ref{lastprop}) 
to produce a perfect, $\Gal(K/\k)$-equivariant, 
skew-Hermitian pairing 
$\spair{\;}{\;}_\cR : \Shmd \times \Shmd \to \Dp$.
By Proposition \ref{linalg} we conclude that 
$\dim_{\rf}(\Sh(A/K)/\Sh(A/K)_\div)[\P_L]$ is even.  This is 
Theorem \ref{pol}.
\end{proof}

\begin{rem}
It is tempting to try to simplify the arguments of this appendix by 
using the pairing of Definition \ref{lastpair} along with the construction 
at the end of the proof of Theorem \ref{flachthm}, to try to produce a perfect, 
skew-symmetric, $G_K$-equivariant pairing $T_p(A) \times T_p(A) \to \Zp(1)$.  
If so, Theorems 1 and 2 of \cite{flach} would give us directly a skew-symmetric 
perfect pairing on $\Shmd$.  Unfortunately, because $\pi^\iota = -\pi$ and 
the different of $\cR/\Zp$ is an odd power of $\p$, one can produce 
in this way (as in the proof of Proposition \ref{sd}) 
a perfect {\em symmetric} pairing, but not a skew-symmetric one.
\end{rem}

\section{The local norm map in the ordinary case}
\label{appb}

In this appendix we  study the cokernel of the local norm map when 
$E$ has ordinary reduction, following and expanding on the proof from
\cite{luro} of some of the results of \cite{mazurav}.  
Our main result is 
Proposition \ref{oni}, which is used to prove Theorem \ref{par4}.

If $\K$ is an algebraic extension of $\Qp$ and $E$ is an elliptic curve 
over $\K$ with good ordinary reduction,
let $E_1(\K)$ denote the kernel of reduction in $E(\K)$, and 
let $U_1(\K)$ denote the units in the ring of integers of $\K$ 
congruent to $1$ modulo the maximal ideal.
We can identify $E_1(\K)$ (resp., $U_1(\K)$) with the maximal ideal 
of $\K$ under the operation given by the formal group of $E$ 
(resp., the formal multiplicative group).  

Suppose now that $\K$ is a finite extension of $\Qp$, with residue 
field $\kappa$.  
Let $u \in \Zp^\times$ be the unit eigenvalue of Frobenius 
acting on the $\ell$-adic Tate module of $E$, for $\ell \ne p$.  
Following \cite{mazurav}, 
we say that $E$ has anomalous reduction if $E(\kappa)[p] \ne 0$, 
or equivalently if $u \equiv 1 \pmod{p}$.  

Fix a totally ramified cyclic extension $\L/\K$ of degree $p^n$.  
Let $\phi$ denote the Frobenius generator of $\Gal(\L^\unr/\L)$; the 
restriction of $\phi$ is the Frobenius generator of 
$\Gal(\K^\unr/\K)$.  
Let $I_{\L/\K} \subset \Z[\Gal(\L/\K)]$ denote the augmentation ideal.

\begin{lem}
\label{bigd}
There is a commutative diagram with exact rows and columns
$$
\xymatrix@C=22pt{
&& 0 \ar[d] & 0 \ar[d] \\
&& \hskip -15pt E_1(\L)/(E_1(\L) \cap I_{\L/\K} U_1(\L^\unr)) \ar^-{N_{\L/\K}}[r]\ar[d] & E_1(\K) \ar[d] \\
0 \ar[r] & \Gal(\L/\K) \ar[r] \ar^{1-u}[d] & U_1(\L^\unr)/I_{\L/\K} U_1(\L^\unr) 
    \ar^-{N_{\L/\K}}[r]\ar^{\phi-u}[d] & U_1(\K^\unr) \ar[r]\ar^{\phi-u}[d] & 0 \\
0 \ar[r] & \Gal(\L/\K) \ar[r] & U_1(\L^\unr)/I_{\L/\K} U_1(\L^\unr) 
    \ar^-{N_{\L/\K}}[r]\ar[d] & U_1(\K^\unr) \ar[r]\ar[d] & 0 \\
&& 0 & 0
}
$$
\end{lem}

\begin{proof}
This is proved on page 239 of \cite{luro}, using an identification 
$$
E_1(\L) \cong \{x \in U_1(\L^\unr) : x^\phi = x^u\}
$$
(see the Lemma on page 237 of \cite{luro}).
\end{proof}

\begin{prop}
\label{e1ind}
Suppose $\K \subset \M \subset \L$ and $[\L:\M] = p$.  Then
$$
\dim_{\Fp}(E_1(\K) / (E_1(\K) \cap N_{\L/\M}E_1(\L))) = 
\begin{cases}
1 & \text{if $E$ has anomalous reduction}, \\
0 & \text{otherwise}.
\end{cases}
$$
\end{prop}

\begin{proof}
Let $G := \Gal(\L/\K)$ and $H := \Gal(\L/\M)$.  
There is a commutative diagram
\begin{equation}
\label{smalld}
\parbox{4in}{
\xymatrix{
E_1(\M)/N_{\L/\M}E_1(\L) \ar^-{\sim}[r] & H/(1-u)H \\
E_1(\K)/N_{\L/\K}E_1(\L) \ar^-{\sim}[r] \ar[u] & G/(1-u)G \ar@{->>}_{\Tr}[u]
}}
\end{equation}
where the horizontal isomorphisms are Corollaries 4.30 and 4.37 of \cite{mazurav}, 
(proved in \cite{luro} by applying the Snake Lemma to the diagram of Lemma
\ref{bigd} for $\L/\M$ and $\L/\K$),
the left-hand vertical map is induced by the 
inclusion of $\K$ into $\M$, and the right-hand vertical map is induced by the 
transfer map $G \to H$.  The commutativity of the diagram follows from 
Lemma \ref{bigd} and the commutativity of 
$$
\xymatrix{
0 \ar[r] & H \ar[r] & U_1(\L^\unr)/I_{\L/\M} U_1(\L^\unr) 
    \ar^-{N_{\L/\M}}[r] & U_1(\M^\unr) \ar[r] & 0 \\
0 \ar[r] & G \ar[r] \ar^{\Tr}[u] & U_1(\L^\unr)/I_{\L/\K} U_1(\L^\unr) 
    \ar^-{N_{\L/\K}}[r]\ar^{N_{\M/\K}}[u] & U_1(\K^\unr) \ar[r]\ar@{^(->}[u] & 0 
}
$$
(see the proof of Lemma 2 of \cite{luro}).

If $E$ has non-anomalous reduction, then $1-u \in \Zp^\times$ 
so the top isomorphism of \eqref{smalld} 
shows that $N_{\L/\M}E_1(\L) = E_1(\M) \supset E_1(\K)$.

If $E$ has anomalous reduction, then $(1-u)H \subset pH = 0$.  
Since $G$ is cyclic, the transfer map is surjective.  
Therefore \eqref{smalld} shows $E_1(\M)/N_{\L/\M}E_1(\L)$ has order $p$, 
and is generated by the image of $E_1(\K)$.  
The proposition follows.
\end{proof}

\begin{prop}
\label{oni}
Suppose that $E$ is defined and has good reduction 
over a subfield $\K^+ \subset \K$ such that $[\K:\K^+] = 2$, 
$\L/\K^+$ is Galois, and $\Gal(\L/\K^+)$ is dihedral.  
If $\K \subset \M \subset \L$ and $[\L:\M] = p$, then 
$$
\dim_{\Fp}(E(\K) / (E(\K) \cap N_{\L/\M}E(\L))) = 
\begin{cases}
2 & \text{if $E$ has anomalous reduction}, \\
0 & \text{otherwise}.
\end{cases}
$$
\end{prop}

\begin{proof}
Let $\kappa$ denote the common residue field of $\K$, $\M$, and $\L$.  
We have a commutative diagram
\begin{equation}
\label{asd}
\parbox{4in}{\xymatrix{
0 \ar[r] & E_1(\L) \ar[r] \ar^{N_{\L/\M}}[d] & E(\L) \ar[r] \ar^{N_{\L/\M}}[d] 
    & E(\kappa) \ar[r] \ar^{p}[d] & 0 \\
0 \ar[r] & E_1(\M) \ar[r] & E(\M) \ar[r] & E(\kappa) \ar[r] & 0
}}
\end{equation}
If $E$ has non-anomalous reduction, then $E(\kappa)$ has order prime to $p$ and 
the proposition follows from Proposition \ref{e1ind}.

Suppose now that $E$ has anomalous reduction.  
Let $H := \Gal(\L/\M)$, and fix $\L^+$ with 
$\K^+ \subset \L^+ \subset \L$, $[\L:\L^+] = 2$.  Let $\M^+ = \M \cap \L^+$.  

Replacing $E/\K^+$ by its quadratic twist by $\K/\K^+$ if necessary, 
we may suppose that $E$ has anomalous reduction over $\K^+$.
We will show that 
\begin{equation}
\label{claim}
\text{$N_{\L/\M} : E_1(\L^+) \to E_1(\M^+)$ is surjective.}
\end{equation}
Assuming this for the moment, choose $x \in E(\K^+)$ 
such that the reduction of $x$ has order $p$ in $E(\kappa)$.  Then 
$\N_{\L/\M}(x) = px \in E_1(\M^+)$ so we 
can find $y \in E_1(\L^+)$ such that $N_{\L/\M}(y) = \N_{\L/\M}(x)$.  
Then $N_{\L/\M}(x-y) = 0$ and the reduction of $x-y$ is nontrivial.  
Therefore, 
since $E(\kappa)[p]$ is cyclic of order $p$, the Snake Lemma applied 
to \eqref{asd} gives an exact sequence
\begin{equation}
\label{anoth}
0 \to E_1(\M)/N_{\L/\M}E_1(\L) \to E(\M)/N_{\L/\M}E(\L) 
    \to E(\kappa)/p E(\kappa) \to 0.
\end{equation}
Using the natural injections 
$E_1(\K)/(E_1(\K) \cap N_{\L/\M}E_1(\L)) \hookto E_1(\M)/N_{\L/\M}E_1(\L)$ 
and $E(\K)/(E(\K) \cap N_{\L/\M}E(\L)) \hookto E(\M)/N_{\L/\M}E(\L)$, 
\eqref{anoth} restricts to an exact sequence
\begin{multline*}
0 \to E_1(\K)/(E_1(\K) \cap N_{\L/\M}E_1(\L)) \\
    \to E(\K)/(E(\K) \cap N_{\L/\M}E(\L)) 
    \to E(\kappa)/p E(\kappa) \to 0.
\end{multline*}
Now the proposition follows from Proposition \ref{e1ind}.

It remains to prove \eqref{claim}.
We consider two cases.  

{\em Case 1: $\K/\K^+$ is unramified.}
Let $v$ be the unit eigenvalue of Frobenius over $\K^+$, so $v^2 = u$.  
Since $E$ has anomalous reduction over $\K^+$, $v \equiv 1 \pmod{p}$.  
Let $\psi$ denote the Frobenius generator of 
$\Gal(\L^\unr/\L^+)$ (note that $(\L^+)^\unr = \L^\unr$), so $\psi^2 = \phi$.
As in Lemma \ref{bigd}, there is a commutative diagram with exact rows and columns
$$
\xymatrix@C=22pt{
&& 0 \ar[d] & 0 \ar[d] \\
&& \hskip -15pt E_1(\L^+)/(E_1(\L^+) \cap I_{\L/\M} U_1(\L^\unr)) 
    \ar^-{N_{\L/\M}}[r]\ar[d] & E_1(\M^+) \ar[d] \\
0 \ar[r] & H \ar[r] \ar^{-1-v}[d] & U_1(\L^\unr)/I_{\L/\M} U_1(\L^\unr) 
    \ar^-{N_{\L/\M}}[r]\ar^{\psi-v}[d] & U_1(\M^\unr) \ar[r]\ar^{\psi-v}[d] & 0 \\
0 \ar[r] & H \ar[r] & U_1(\L^\unr)/I_{\L/\M} U_1(\L^\unr) 
    \ar^-{N_{\L/\M}}[r]\ar[d] & U_1(\M^\unr) \ar[r]\ar[d] & 0 \\
&& 0 & 0
}
$$
The proof is the same as the proof in \cite{luro} of Lemma \ref{bigd}.  
The only point to notice is the map $-1-v$ on the left, which arises because 
if $\pi$ is a uniformizing parameter of $\L^+$ and $h \in H$, then 
$\psi h \psi^{-1} = h^{-1}$ on $\L$, so 
$$
(\pi^{h-1})^{1+\psi} = \pi^{h - 1 + \psi h^{-1}-\psi} = \pi^{h+h^{-1}-2} 
    = (\pi^{h-1})^{1-h^{-1}} \in I_{\L/\M}U^1(\L).
$$
Since the left-most horizontal maps send $h \mapsto \pi^{h-1}$, this shows that the 
left-hand square commutes (see \cite{luro} page 239).
Since $p \ne 2$, $-1-v \in \Zp^\times$, and \eqref{claim} now follows from the Snake Lemma 
in this case.  

{\em Case 2: $\K/\K^+$ is ramified.}
In this case $\L^\unr/(\L^+)^\unr$ is a quadratic extension.  
Taking $\Gal(\L^\unr/(\L^+)^\unr)$-invariants in the diagram of Lemma 
\ref{bigd} (applied to $\L/\M$) 
gives a new diagram with exact rows and columns.  
The top row of the new diagram is 
$$
E_1(\L^+)/(E_1(\L^+) \cap I_{\L/\M} U_1((\L^+)^\unr)) \map{N_{\L/\M}} E_1(\M^+),
$$
and the left-hand column is $0 \to 0$
since $\Gal(\L^\unr/(\L^+)^\unr)$ acts on $H$ by $-1$.  
Now the Snake Lemma applied to 
this new diagram proves \eqref{claim} in this case.  
This completes the proof of the proposition.
\end{proof}


\begin{thebibliography}{MR1}

\bibitem[BD]{bd}
   M.\ Bertolini and H.\ Darmon, 
   Non-triviality of families of Heegner points and ranks of Selmer groups 
   over anticyclotomic towers. {\em J.\ Ramanujan Math.\ Soc.}
   {\bf 13} (1998) 15--24.

\bibitem[BK]{bk} 
   Bloch, S., Kato, K.: $L$-functions and Tamagawa numbers of
   motives.  In: The 
   Grothendieck Festschrift (Vol.\ I), P.\ Cartier, et al., eds., 
   {\em Prog.\ in Math.}\ {\bf 86}, Birkh\"auser, Boston (1990) 333--400.
   
\bibitem[Br]{brown}
    K.\ Brown, Cohomology of groups, {\em Grad.\ Texts in Math.} {\bf 87},
    Springer, New York (1982).

\bibitem[Ca1]{cassels}
   J.W.S.\ Cassels, Arithmetic on curves of genus $1$. IV. Proof of the Hauptvermutung. 
   {\em J.\ Reine Angew.\ Math.} {\bf 211} (1962) 95--112.

\bibitem[Ca2]{casexp}
   J.W.S.\ Cassels, Diophantine equations with special reference to elliptic curves. 
   {\em J.\ London Math.\ Soc.} {\bf 41} (1966) 193--291.

\bibitem[Co]{cornut}
   C.\ Cornut, Mazur's conjecture on higher Heegner points.
   {\em Invent.\ Math.} {\bf 148} (2002) 495--523, MR1908058, Zbl pre01777253.

\bibitem[F]{flach}
   M.\ Flach, A generalisation of the Cassels-Tate pairing. 
   {\em J.\ Reine\ Angew.\ Math.} {\bf 412} (1990) 113--127. 

\bibitem[Hb]{howard}
   B.\ Howard, The Heegner point Kolyvagin system. 
   {\em Compositio Math.} {\bf 140} (2004) 1439--1472.
   
\bibitem[He]{howe}
   E.\ Howe, Isogeny classes of abelian varieties with no principal polarizations. 
   In: Moduli of abelian varieties (Texel Island, 1999), {\em Progr.\ Math.}
   {\bf 195}, Birkh\"auser, Basel (2001) 203--216.

\bibitem[K]{kim}
   B-D.\ Kim, The parity conjecture and algebraic functional 
   equations for elliptic curves at supersingular 
   reduction primes.  Thesis, Stanford University (2005).

\bibitem[LR]{luro}
   J.\ Lubin and M.\ Rosen, The norm map for ordinary abelian varieties. 
   {\em J.\ Algebra} {\bf 52} (1978) 236--240.

\bibitem[M]{mazurav}
   B.\ Mazur, Rational points of abelian varieties with values in towers 
   of number fields. {\em Invent.\ Math.} {\bf 18} (1972) 183--266.

\bibitem[MR1]{kolysys}
   B.\ Mazur, K.\ Rubin, Kolyvagin systems.  
   {\em Mem.\ Amer.\ Math.\ Soc.} {\bf 168} (2004) no.~799.

\bibitem[MR2]{organization}
   B.\ Mazur, K.\ Rubin, Organizing the arithmetic of elliptic curves.  
   {\em Advances in Math.} {\bf 198} (2005) 504--546.

\bibitem[MR3]{finding}
   B.\ Mazur, K.\ Rubin, Finding large Selmer groups.  
   {\em J.\ Differential Geom.} {\bf 70} (2005) 1--22.

\bibitem[MRS]{prim}
   B.\ Mazur, K.\ Rubin, A.\ Silverberg, Twisting commutative algebraic groups.  
   Preprint available at 
   {\tt http://arxiv.org/math/0609066}

\bibitem[Mi]{milne}
     J.\ S.\ Milne, {\em On the arithmetic of abelian varieties}. 
     Invent.\ Math.\ {\bf 17} (1972) 177--190.     

\bibitem[N1]{nekovarpc}
   J.\ Nekov\'a\u{r}, On the parity of ranks of Selmer groups. II,
   {\em  C.\ R.\ Acad.\ Sci.\ Paris S\'er.\ I Math.} {\bf 332} (2001) 99--104. 
   
\bibitem[N2]{nekovar}
   J.\ Nekov\'a\u{r}, Selmer complexes (Fourth version, May 2006).  
   Preprint available at 
   {\tt  http://www.math.jussieu.fr/$\sim$nekovar/pu/}.

\bibitem[R]{EulerSystems}
   K.\ Rubin, Euler Systems.  {\em Annals of Math. Studies} {\bf 147},
   Princeton University Press, Princeton (2000).

\bibitem[S]{stein}
   W.\ Stein, 
   Shafarevich-Tate groups of nonsquare order. In: Modular curves and abelian 
   varieties, {\em Progr.\ Math.} {\bf 224}, Birkh\"auser, Basel (2004) 277--289.

\bibitem[T1]{tate}
    J.\ Tate, WC-groups over $\p$-adic fields.  {\em S\'eminaire Bourbaki} 
    {\bf 4}, Exp.\ 156, Paris: Soc.\ Math.\ France (1995) 265--277.

\bibitem[T2]{tate2}
    J.\ Tate, Duality theorems in Galois cohomology over number fields.  
    In: Proc.\ Internat.\ Congr.\ Mathematicians (Stockholm, 1962)
    Inst.\ Mittag-Leffler, Djursholm (1963) 288--295.

\bibitem[V]{vatsal}
   V.\ Vatsal, Uniform distribution of Heegner points. 
   {\em Invent.\ Math.} {\bf 148} (2002) 1--46. 

\bibitem[W]{weil}
    A.\ Weil,
    Adeles and algebraic groups. {\em Progress in Math.} {\bf 23},
    Birkh\"auser, Boston (1982).
    
\end{thebibliography}
\end{document}